\documentclass[article]{amsart}

\usepackage{amsmath,amsthm,amssymb,amsfonts,amscd}
\usepackage{url}
\usepackage{color}
\usepackage{setspace}
\usepackage{enumerate}
\usepackage{enumitem}
\usepackage{mathrsfs}

\usepackage{epsfig,graphicx,pinlabel}

\definecolor{mygreen}{rgb}{0.0,0.5,0.1}

\theoremstyle{plain}
\newtheorem{theorem}{Theorem}
\newtheorem{proposition}[theorem]{Proposition}
\newtheorem{condition}[theorem]{Condition}

\theoremstyle{definition}
\newtheorem{definition}[theorem]{Definition}
\newtheorem{example}[theorem]{Example}
\newtheorem{examples}[theorem]{Examples}
\newtheorem{remark}[theorem]{Remark}
\newtheorem{problem}{Problem}

\numberwithin{theorem}{section}
\numberwithin{equation}{section}

\newcommand{\N}{\mathbb{N}}
\newcommand{\Z}{\mathbb{Z}}

\newcommand{\R}{\mathbb{R}}
\newcommand{\C}{\mathbb{C}}
\newcommand{\HH}{\mathbb{H}}
\newcommand{\PP}{\mathbb{P}}

\newcommand{\SL}{\mathrm{SL}}
\newcommand{\GL}{\mathrm{GL}}
\newcommand{\SO}{\mathrm{SO}}
\newcommand{\OO}{\mathrm{O}}
\newcommand{\PO}{\mathrm{PO}}
\newcommand{\PSL}{\mathrm{PSL}}
\newcommand{\PGL}{\mathrm{PGL}}

\newcommand{\g}{\mathfrak{g}}

\newcommand{\AdS}{\mathrm{AdS}}

\newcommand{\Hom}{\mathrm{Hom}}
\newcommand{\ie}{i.e.\ }
\newcommand{\eg}{e.g.\ }
\newcommand{\resp}{resp.\ }

\newcommand{\CC}{\mathcal{C}}

\title[Geometric structures and representations of discrete groups]{Geometric structures and representations of discrete groups}

\author[F. Kassel]{Fanny Kassel}

\thanks{Partially supported by the European Research Council under Starting Grant 715982 (DiGGeS)}

\address{CNRS and Institut des Hautes \'Etudes Scientifiques, Laboratoire Alexander Grothendieck, 35 route de Chartres, 91440 Bures-sur-Yvette, France}

\email{kassel@ihes.fr}

\begin{document}

\maketitle

\vspace{-1cm}

\begin{abstract}
We describe recent links between two topics: geometric structures on manifolds in the sense of Ehresmann and Thurston, and dynamics ``at infinity'' for representations of discrete groups into Lie groups.
\end{abstract}

\section{Introduction}

The goal of this survey is to report on recent results relating geometric structures on manifolds to dynamical aspects of representations of discrete groups into Lie groups, thus linking geometric topology to group theory and dynamics.

\subsection{Geometric structures} \label{subsec:intro-GXstruct}

\begin{figure}[ht!]
\centering
\includegraphics[scale=0.6]{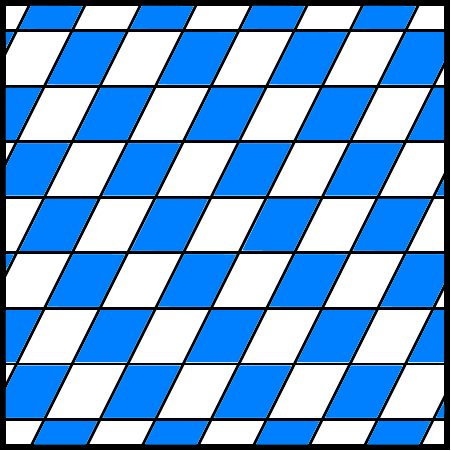}\hspace{1cm}\includegraphics[scale=0.6]{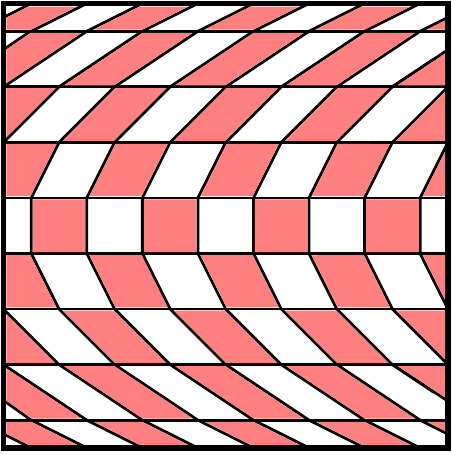}\hspace{1cm}\includegraphics[scale=0.195]{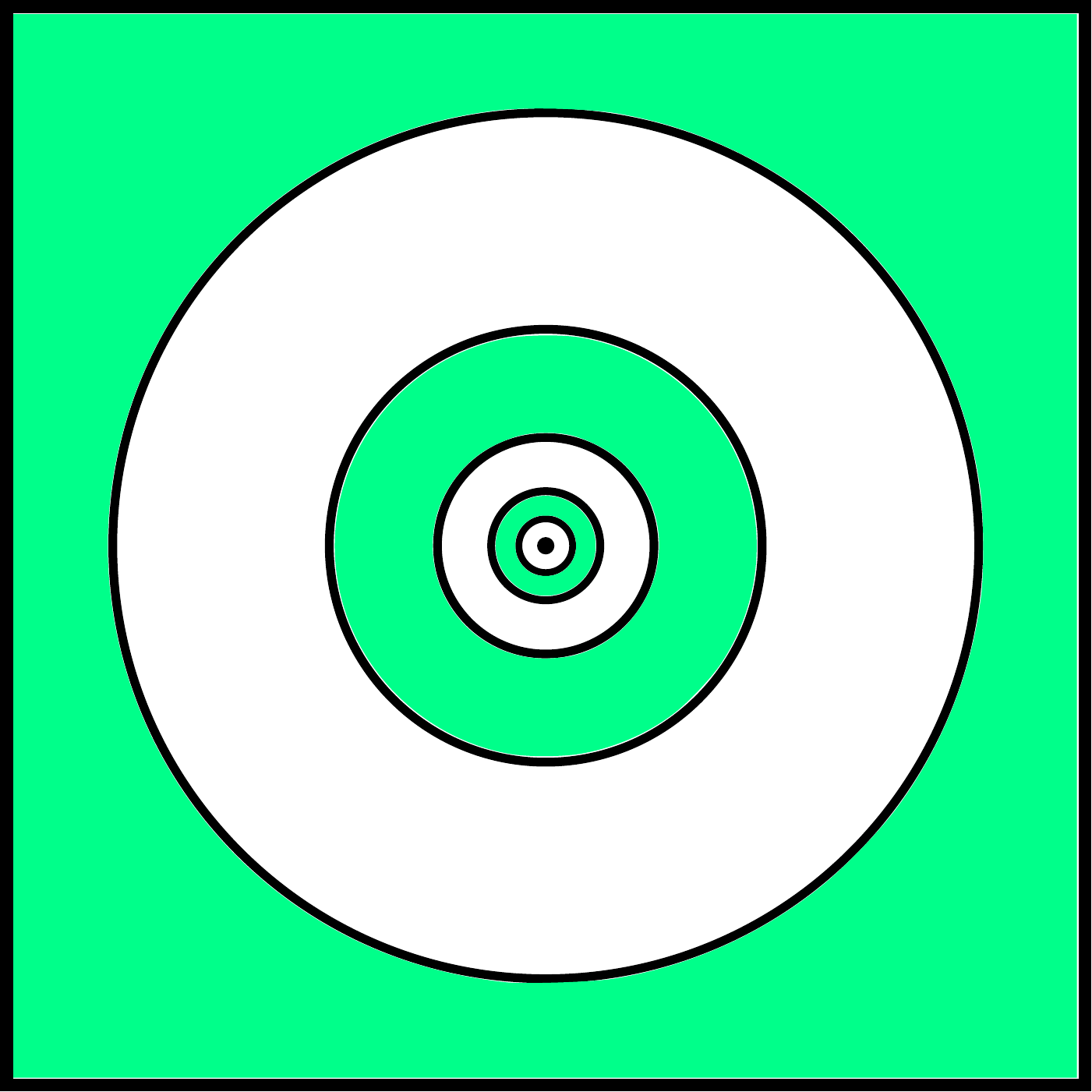}
\vspace{-0.1cm}
\caption{Tilings of $X=\R^2$ showing the three $\Gamma$-actions in Example~\ref{ex:torus}}
\vspace{-0.4cm}
\label{fig-PavAff}
\end{figure}

The first topic of this survey is geometric structures on manifolds.
Here is a concrete example as illustration (see Figure~\ref{fig-PavAff}).

\begin{example} \label{ex:torus}
Consider a two-dimensional torus~$T$.

(1) We can view $T$ as the quotient of the Euclidean plane $X=\R^2$ by $\Gamma=\Z^2$, which is a discrete subgroup of the isometry group $G\!=\!\OO(2)\ltimes\R^2$\,of\,$X$ (acting by linear isometries and translations).
Viewing $T$ this way provides it with a Riemannian metric and a notion of parallel lines, length, angles, etc.
We say $T$ is endowed with a \emph{Euclidean} (or \emph{flat}) \emph{structure}, or a $(G,X)$-structure with $(G,X) = (\OO(2)\ltimes\R^2,\R^2)$.

(2) Here is a slightly more involved way to view $T$: we can see it as the quotient of the affine plane $X=\R^2$ by the group $\Gamma$ generated by the translation of vector $(\begin{smallmatrix} 1\\ 0\end{smallmatrix})$ and the affine transformation with linear part $(\begin{smallmatrix} 1 & 1\\ 0 & 1\end{smallmatrix})$ and translational part $(\begin{smallmatrix} 0\\ 1\end{smallmatrix})$.
This group $\Gamma$ is now a discrete subgroup of the affine group $G=\GL(2,\R)\ltimes\R^2$.
Viewing $T$ this way still provides it with a notion of parallel lines and even of geodesic, but no longer with a notion of length or angle or speed of geodesic.
We say $T$ is endowed with an \emph{affine structure}, or a $(G,X)$-structure with $(G,X) = (\GL(2,\R)\ltimes\R^2,\R^2)$.

(3) There are many ways to endow $T$ with an affine structure.
Here is a different one: we can view $T$ as the quotient of the open subset $\mathcal{U}=\R^2\smallsetminus\{0\}$ of $X=\R^2$ by the discrete subgroup $\Gamma$ of $G=\GL(2,\R)\ltimes\R^2$ generated by the homothety $(\begin{smallmatrix} 2 & 0\\ 0 & 1/2\end{smallmatrix})$.
This still makes $T$ ``locally look like'' $X=\R^2$, but now the image in~$T$ of an affine geodesic of~$X$ pointing towards the origin is \emph{incomplete} (it circles around in~$T$ with shorter and shorter period and disappears in a finite amount of time).
\end{example}

As in Example~\ref{ex:torus}, a key idea underlying a large part of modern geometry is the existence of \emph{model geometries} which various manifolds may locally carry.
By definition, a model geometry is a pair $(G,X)$ where $X$ is a manifold (\emph{model space}) and $G$ a Lie group acting transitively on~$X$ (\emph{group of symmetries}).
In Example~\ref{ex:torus} we encountered $(G,X) = (\OO(n)\ltimes\R^n,\R^n)$ and $(G,X) = (\GL(n,\R)\ltimes\R^n,\R^n)$, corresponding respectively to \emph{Euclidean geometry} and \emph{affine geometry}.
Another important example is $X=\HH^n$ (the $n$-dimensional real hyperbolic space) and $G=\PO(n,1)=\OO(n,1)/\{\pm I\}$ (its group of isometries), corresponding to \emph{hyperbolic geometry}.
(For $n=2$ we can see $X$ as the upper half-plane and $G$, up to index two, as $\PSL(2,\R)$ acting by homographies.)
We refer to Table~\ref{table1} for more~examples.

\begin{figure}[ht!]
\labellist
\small\hair 2pt
\pinlabel \textcolor{red}{$\scriptstyle g\in G$} at 460 135
\pinlabel $X$ at 560 160
\pinlabel $M$ at 25 145
\endlabellist
\centering
\includegraphics[scale=0.35]{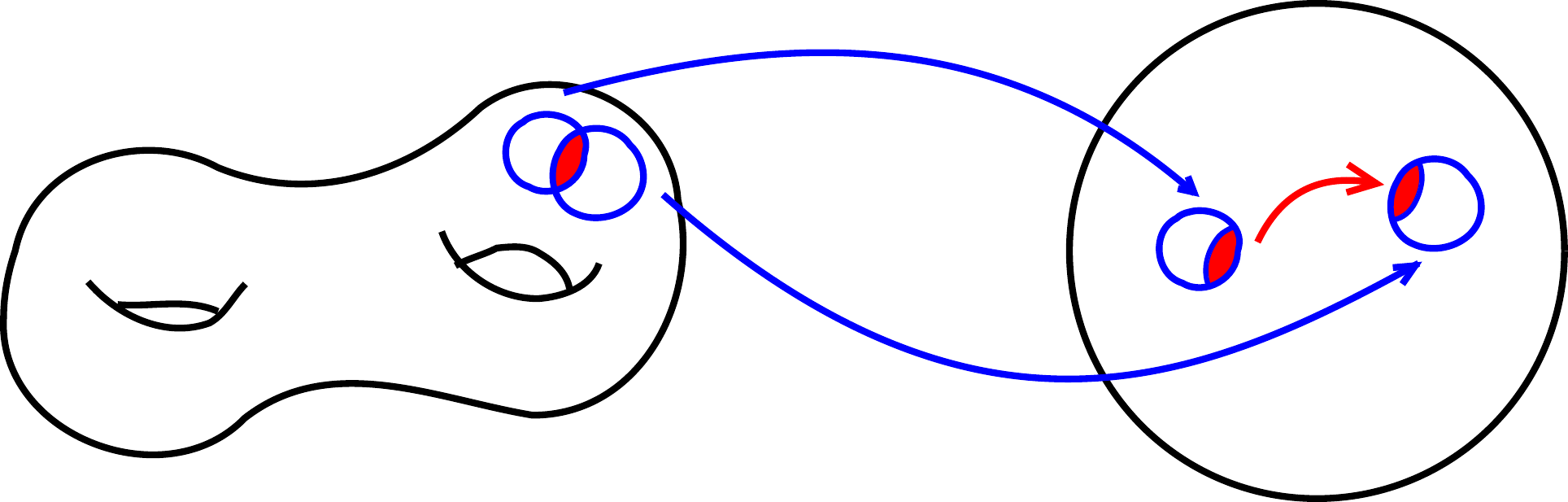}
\vspace{-0.5cm}
\caption{Charts defining a $(G,X)$-structure on~$M$}
\vspace{-0.2cm}
\label{fig-GXstruct}
\end{figure}

The idea that a manifold~$M$ locally carries the geometry $(G,X)$ is formalized by the notion of a \emph{$(G,X)$-structure} on~$M$: by definition, this is a maximal atlas of coordinate charts on~$M$ with values in~$X$ such that the transition maps are given by elements of~$G$ (see Figure~\ref{fig-GXstruct}).
Note that this is quite similar to a manifold structure on~$M$, but we now require the charts to take values in~$X$ rather than~$\R^n$, and the transition maps to be given by elements of~$G$ rather than diffeomorphisms of~$\R^n$.
\vspace{-0.2cm}
Although general $(G,X)$-structures may display pathological behavior (see \cite{gol-book}), in this survey we will restrict to the two ``simple'' types of $(G,X)$-structures appearing in Example~\ref{ex:torus}, to which we shall give names to facilitate the discussion:
\begin{itemize}[leftmargin=0.8cm]
  \item \textbf{Type C (``complete'')}: $(G,X)$-structures that identify $M$ with a quotient of~$X$ by a discrete subgroup $\Gamma$ of~$G$ acting properly discontinuously;
  \item \textbf{Type U (``incomplete but still uniformizable'')}: $(G,X)$-structures that identify $M$ with a quotient of some proper open subset $\mathcal{U}$ of~$X$ by a discrete subgroup $\Gamma$ of~$G$ acting properly discontinuously.
\end{itemize}
Setting $\mathcal{V}=X$ or~$\mathcal{U}$ as appropriate, we then have coverings $\widetilde{M}\simeq\widetilde{\mathcal{V}}\to\mathcal{V}\to\nolinebreak\Gamma\backslash\mathcal{V}\simeq M$ (where\; $\widetilde{ }$\; denotes universal covers).
The charts on~$M$ are obtained by taking preimages in $\mathcal{V}\subset X$ of open subsets of~$M$.
Moreover, the basic theory of covering groups gives a natural group homomorphism $\mathrm{hol} : \pi_1(M)\to G$ with image~$\Gamma$ and kernel $\pi_1(\mathcal{V})$, called the \emph{holonomy}.

In this survey, we use the phrase \emph{geometric structures} for $(G,X)$-structures.
We shall not detail the rich historical aspects of geometric structures here; instead, we refer to the excellent surveys \cite{gol-icm10,gol-historical,gol-book}.
We just mention that the notion of model geometry has its origins in ideas of Lie and Klein, formulated in Klein's 1872 Erlangen program.
Influenced by these ideas and those of Poincar\'e, Cartan and others, Ehresmann \cite{ehr37} initiated a general study of geometric structures in 1935.
Later, geometric structures were greatly promoted by Thurston's revolutionary work \cite{thu80}.

\subsection{Classifying geometric structures}

The fundamental problem in the theory of geometric structures is their classification, namely:

\begin{problem} \label{problemA}
Given a manifold~$M$,
\begin{enumerate}
  \item Describe which model geometries $(G,X)$ the manifold~$M$ may locally carry;
  \item For a fixed model $(G,X)$, describe all possible $(G,X)$-structures on~$M$.
\end{enumerate}
\end{problem}

We refer to \cite{gol-icm10} for a detailed survey of Problem~\ref{problemA} with a focus on dimensions two and three, and to \cite{ky05} for a special case.

Problem~\ref{problemA}.(1) asks how the global topology of~$M$ determines the geometries that it may locally carry.
This has been the object of deep results, among which:
\begin{itemize}[leftmargin=0.8cm]
  \item the classical \emph{uniformization theorem}: a closed Riemann surface~may carry a Euclidean, a spherical, or a hyperbolic structure, depending on its genus;
  \item Thurston's \emph{hyperbolization theorem}: a large class of $3$-dimensional manifolds, defined in purely topological terms, may carry a hyperbolic structure;
  \item more generally, Thurston's \emph{geometrization program} (now Perelman's theorem): any closed orientable $3$-dimensional manifold may be decomposed into pieces, each admitting one of eight model geometries (see \cite{bon02}).
\end{itemize}

Problem~\ref{problemA}.(2) asks to describe the \emph{deformation space} of $(G,X)$-structures on~$M$.
In the simple setting of Example~\ref{ex:torus}, this space is already quite rich (see \cite{ny73}).
For hyperbolic structures on a closed Riemann surface of genus $\geq 2$ (Example~\ref{ex:holonomy-Teich}), Problem~\ref{problemA}.(2) gives rise to the fundamental and wide-ranging \emph{Teichm\"uller theory}.

\subsection{Representations of discrete groups} \label{subsec:intro-holonomy}

The second topic of this survey is representations (\ie group homomorphisms) of discrete groups (\ie countable groups) to Lie groups~$G$, and their dynamics ``at infinity''.
We again start with an example.

\begin{figure}
\centering
\vspace{-0.2cm}
\includegraphics[scale=0.28]{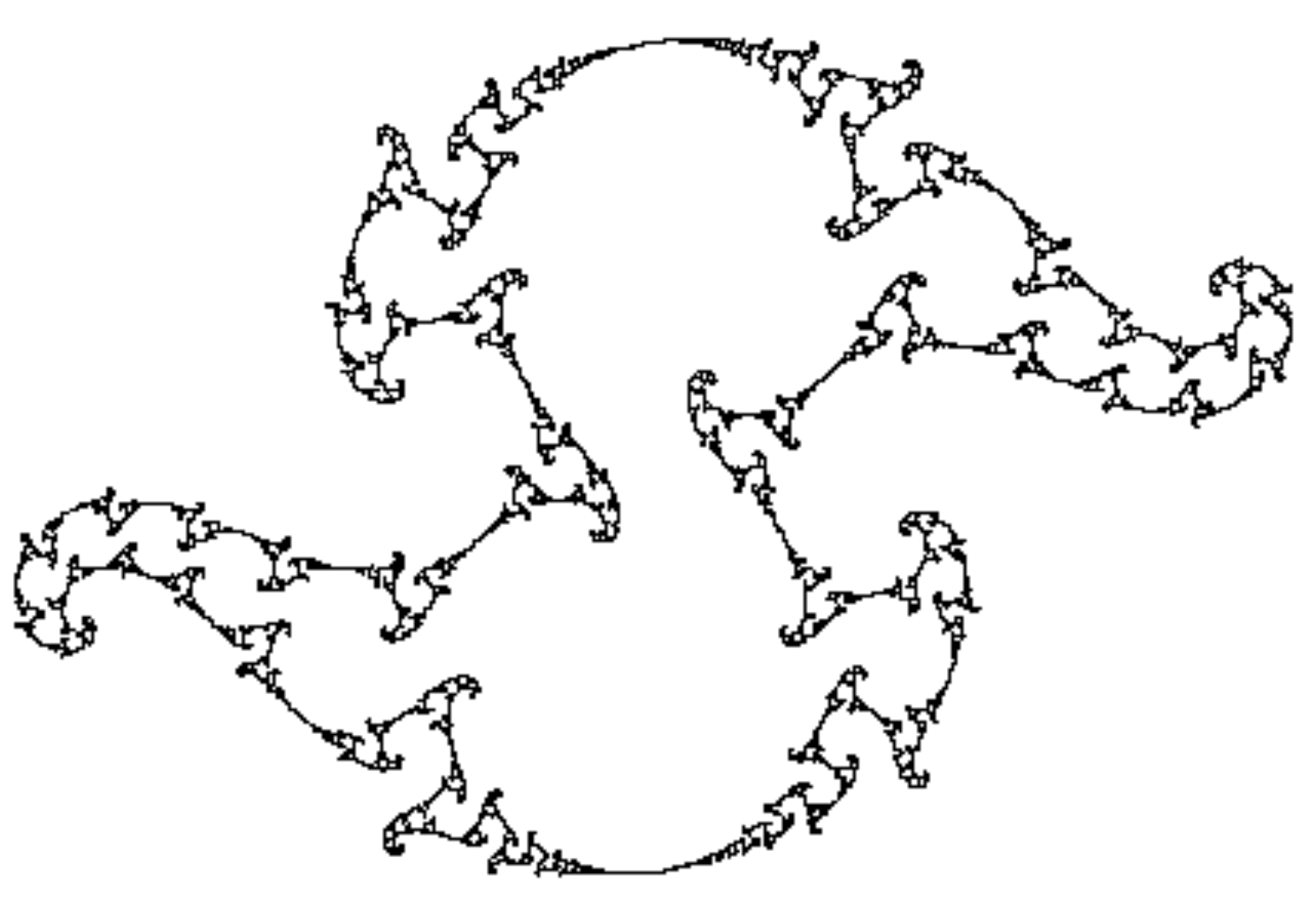}
\vspace{-0.4cm}
\caption{The limit set of a quasi-Fuchsian group in $\partial_{\infty}\HH^3 \simeq \C \cup \{\infty\}$}
\vspace{-0.3cm}
\label{fig-qF}
\end{figure}

\begin{example} \label{ex:qF-intro}
Let $\Gamma = \pi_1(S)$ where $S$ is a closed oriented Riemann surface of genus $\geq 2$.
By the uniformization theorem, $S$ carries a complete (``type~C'') hyperbolic structure, which yields a holonomy representation $\Gamma\to\PSL(2,\R)$ as in Section~\ref{subsec:intro-GXstruct}.
Embedding $\PSL(2,\R)$ into $G=\PSL(2,\C)$, we obtain a representation $\rho : \Gamma\to G$, called \emph{Fuchsian}, and an associated action of $\Gamma$ on the hyperbolic space $X=\HH^3$ and on its boundary at infinity $\partial_{\infty}\HH^3 = \widehat{\C}$ (the Riemann sphere).
The \emph{limit set}~of~$\rho(\Gamma)$ in $\widehat{\C}$ is the set of accumulation points of $\rho(\Gamma)$-orbits of~$X$; it is a circle in the~sphere $\widehat{\C}$.
Deforming $\rho$ slightly yields a new representation $\rho' : \Gamma\to G$,~called \emph{quasi-Fuchsian}, which is still faithful, with discrete image, and whose limit set in $\widehat{\C}$ is still a topological circle (now ``wiggly'', see Figure~\ref{fig-qF}).
The action of $\rho'(\Gamma)$ is chaotic on the limit set (\eg all orbits are dense) and properly discontinuous on~its~complement.
\end{example}

\vspace{0.2cm}

Example~\ref{ex:qF-intro} plays a central role in the theory of Kleinian groups and in Thurston's geometrization program; it was extensively studied by Ahlfors, Beardon, Bers, Marden, Maskit, Minsky, Sullivan, Thurston, and many others.

In this survey we report on various generalizations of Example~\ref{ex:qF-intro}, for representations of discrete groups $\Gamma$ into semisimple Lie groups~$G$ which are faithful (or with finite kernel) and whose images are discrete subgroups of~$G$.
While in Example~\ref{ex:qF-intro} the group $G=\PSL(2,\C)$ has real rank one (meaning that its Riemannian symmetric space $\HH^3$ has no flat region beyond geodesics), we also wish to consider the case that $G$ has \emph{higher real rank}, \eg $\PGL(d,\R)$ with $d\geq 3$.
In general, semisimple groups $G$ tend to have very different behavior depending on whether their real rank is one or higher; for instance, the \emph{lattices} of~$G$ (\ie the discrete subgroups of finite covolume for the Haar measure) may display some forms of flexibility in real rank one, but exhibit strong rigidity phenomena in higher real rank.
Beyond lattices, the landscape of discrete subgroups of~$G$ is somewhat understood in real rank one (at least several important classes of discrete subgroups have been identified for their good geometric, topological, and dynamical properties, see Section~\ref{subsec:cc-rank1}), but it remains very mysterious in higher real rank.
We shall explain some recent attempts at understanding it better.

One interesting aspect is that, even when $G$ has higher real rank, discrete subgroups of~$G$ of infinite covolume may be nonrigid and in fact admit large deformation spaces.
In particular, as part of \emph{higher Teichm\"uller theory}, there has recently been an active and successful effort to find large deformation spaces of faithful and discrete representations of surface groups $\pi_1(S)$ into higher-rank semisimple~$G$ which share some of the rich features of the Teichm\"uller space of~$S$ (see Sections \ref{subsec:Higher-Teich} and~\ref{sec:geom-for-Ano}, and \cite{biw14,wie-icm}).
Such features also include dynamics ``at infinity'' as in Example~\ref{ex:qF-intro}, which are encompassed by a notion of \emph{Anosov representation} \cite{lab06} (see Section~\ref{sec:Anosov}).

\subsection{Flag varieties and boundary maps} \label{subsec:intro-G/P}

Let us be a bit more precise.
Given a representation $\rho : \Gamma\to G$, by dynamics `at infinity'' we mean the dynamics of the action of~$\Gamma$ via~$\rho$ on some \emph{flag varieties} $G/P$ (where $P$ is a parabolic subgroup), seen as ``boundaries'' of~$G$ or of its Riemannian symmetric space $G/K$.
In Example~\ref{ex:qF-intro} we considered a rank-one situation where $G=\PSL(2,\C)$ and $G/P = \partial_{\infty}\HH^3 = \widehat{\C}$.
A typical higher-rank situation that we have in mind is $G=\PGL(d,\R)$ with $d\geq 3$ and $G/P=\mathrm{Gr}_i(\R^d)$ (the Grassmannian of $i$-planes in~$\R^d$) for some $1\leq i\leq d-1$.

In the work of Mostow, Margulis, Furstenberg, and others, rigidity results have often relied on the construction of $\Gamma$-equivariant measurable maps from or to $G/P$.
More recently, in the context of higher Teichm\"uller theory \cite{biw10,fg06,lab06}, it has proved important to study continuous equivariant \emph{boundary maps} which embed the boundary $\partial_{\infty}\Gamma$ of a Gromov hyperbolic group $\Gamma$ (\ie the visual boundary of the Cayley graph of~$\Gamma$) into $G/P$.
Such boundary maps $\xi : \partial_{\infty}\Gamma\to G/P$ define a closed invariant subset $\xi(\partial_{\infty}\Gamma)$ of $G/P$, the \emph{limit set}, on which the dynamics of the action by $\Gamma$ accurately reflects the intrinsic chaotic dynamics of $\Gamma$ on $\partial_{\infty}\Gamma$.
These boundary maps may be used to transfer the Anosov property of the intrinsic geodesic flow of~$\Gamma$ into some uniform contraction/expansion properties for a flow on a natural flat bundle associated to $\rho$ and $G/P$ (see Section~\ref{sec:Anosov}).
They may also define some open subsets $\mathcal{U}$ of $G/P$ on which the action of $\Gamma$ is properly discontinuous, by removing an ``extended limit set'' $\mathscr{L}_{\rho(\Gamma)} \supset \xi(\partial_{\infty}\Gamma)$ (see Sections \ref{sec:dyn-properties-hol}, \ref{sec:geom-for-Ano}, \ref{sec:proj-cc}); this generalizes the domains of discontinuity in the Riemann sphere of Example~\ref{ex:qF-intro}.

For finitely generated groups $\Gamma$ that are not Gromov hyperbolic, one can still define a boundary $\partial_{\infty}\Gamma$ in several natural settings, \eg as the visual boundary of some geodesic metric space on which $\Gamma$ acts geometrically, and the approach considered in this survey can then be summarized by the following general problem.

\begin{problem} \label{problemB}
Given a discrete group~$\Gamma$ with a boundary $\partial_{\infty}\Gamma$, and a Lie group~$G$ with a boundary $G/P$, identify large (\eg open in $\Hom(\Gamma,G)$) classes of faithful and discrete representations $\rho : \Gamma\to G$ for which there exist continuous $\rho$-equivariant boundary maps $\xi : \partial_{\infty}\Gamma\to G/P$.
Describe the dynamics of $\Gamma$ on $G/P$ via~$\rho$.
\end{problem}

\subsection{Goal of the paper}

We survey recent results on $(G,X)$-structures (Problem~\ref{problemA}) and on representations of discrete groups (Problem~\ref{problemB}), making links between the two topics.
In one direction, we observe that various types of $(G,X)$-structures have holonomy representations that are interesting for Problem~\ref{problemB}.
In the other direction, starting with representations that are interesting for Problem~\ref{problemB} (Anosov representations), we survey recent constructions of associated $(G,X)$-structures.
These results tend to indicate some deep interactions between the geometry of $(G,X)$-manifolds and the dynamics of their holonomy representations, which largely remain to be explored.
We hope that they will continue to stimulate the development of rich theories in the future.

\subsection*{Organization of the paper}

In Section~\ref{sec:holonomy} we briefly review the notion of a holonomy representation.
In Section~\ref{sec:dyn-properties-hol} we describe three important families of $(G,X)$-structures for which boundary maps into flag varieties naturally appear.
In Section~\ref{sec:Anosov} we define Anosov representations and give examples and characterizations.
In Section~\ref{sec:geom-for-Ano} we summarize recent constructions of geometric structures associated to Anosov representations.
In Section~\ref{sec:proj-cc} we discuss a situation in which the links between geometric structures and Anosov representations are particularly tight, in the context of convex projective geometry.
In Section~\ref{sec:hol-are-Ano} we examine an instance of $(G,X)$-structures for a nonreductive Lie group~$G$, corresponding to affine manifolds and giving rise to affine Anosov representations.
We conclude with a few remarks.

\subsection*{Acknowledgements}

I would like to heartily thank all the mathematicians who helped, encouraged, and inspired me in the past ten years; the list is too long to include here.
I am very grateful to all my coauthors, in particular those involved in the work discussed below: Jeffrey Danciger (\S \ref{sec:geom-for-Ano},\,\ref{sec:proj-cc},\,\ref{sec:hol-are-Ano}), Fran\c{c}ois Gu\'eritaud (\S \ref{sec:Anosov},\,\ref{sec:geom-for-Ano},\,\ref{sec:proj-cc},\,\ref{sec:hol-are-Ano}), Olivier Guichard (\S \ref{sec:Anosov},\,\ref{sec:geom-for-Ano}), Rafael~Potrie~(\S \ref{sec:Anosov}), and Anna Wienhard (\S \ref{sec:Anosov},\,\ref{sec:geom-for-Ano}).
I warmly thank J.-P.~Burelle, J.~Danciger, O.~Guichard, and S.~Maloni for reading earlier versions of this text and making many valuable comments and suggestions, and R.~Canary and W.~Goldman for kindly answering my questions.

\section{Holonomy representations} \label{sec:holonomy}

Let $G$ be a real Lie group acting transitively, faithfully, analytically on a manifold~$X$, as in Table~\ref{table1}.
In Section~\ref{subsec:intro-GXstruct} we defined holonomy representations for certain types of $(G,X)$-structures.
We now give a short review of the notion~in~general.
\begin{table}[h!]
\centering
\begin{tabular}{|p{3cm}|p{0.8cm}|p{3.8cm}|p{3.9cm}|}
\hline
\centering Type of geometry & \centering $X$ & \centering $G$ & \centering $H$
\tabularnewline
\hline
\centering Real projective & \centering $\PP^n(\R)$ & \centering $\PGL(n+1,\R)$ & stab.\,in $G$ of a line of $\R^{n+1}$\!\!
\tabularnewline
\centering Affine & \centering $\R^n$ & \centering $\mathrm{Aff}(\R^n)=\GL(n,\R)\ltimes\R^n$ & \centering $\GL(n,\R)$
\tabularnewline
\centering Euclidean & \centering $\R^n$ & \centering $\mathrm{Isom}(\R^n)=\OO(n)\ltimes\R^n$ & \centering $\OO(n)$
\tabularnewline
\centering Real hyperbolic & \centering $\HH^n$ & \centering $\mathrm{Isom}(\HH^n)=\PO(n,1)$ & \centering $\OO(n)$
\tabularnewline
\centering Spherical & \centering $\mathbb{S}^n$ & \centering $\mathrm{Isom}(\mathbb{S}^n)=\OO(n+1)$ & \centering $\OO(n)$
\tabularnewline
\centering Complex projective & \centering $\PP^n(\C)$ & \centering $\PGL(n+1,\C)$ & stab.\,in $G$ of a line of $\C^{n+1}$\!\!
\tabularnewline
\hline
\end{tabular}
\medskip
\caption{Some examples of model geometries $(G,X)$, where $X\simeq G/H$}
\vspace{-0.5cm}
\label{table1}
\end{table}

\begin{figure}[ht!]
\labellist
\small\hair 2pt
\pinlabel $M$ at 25 145
\pinlabel $\scriptstyle m$ at 100 80
\pinlabel \textcolor{blue}{$\mathcal{U}$} at 112 110
\pinlabel \textcolor{blue}{$\varphi$} at 280 164
\pinlabel $X$ at 560 160
\pinlabel \textcolor{red}{$\scriptstyle g\in G$} at 450 30
\endlabellist
\centering
\vspace{-0.4cm}
\includegraphics[scale=0.35]{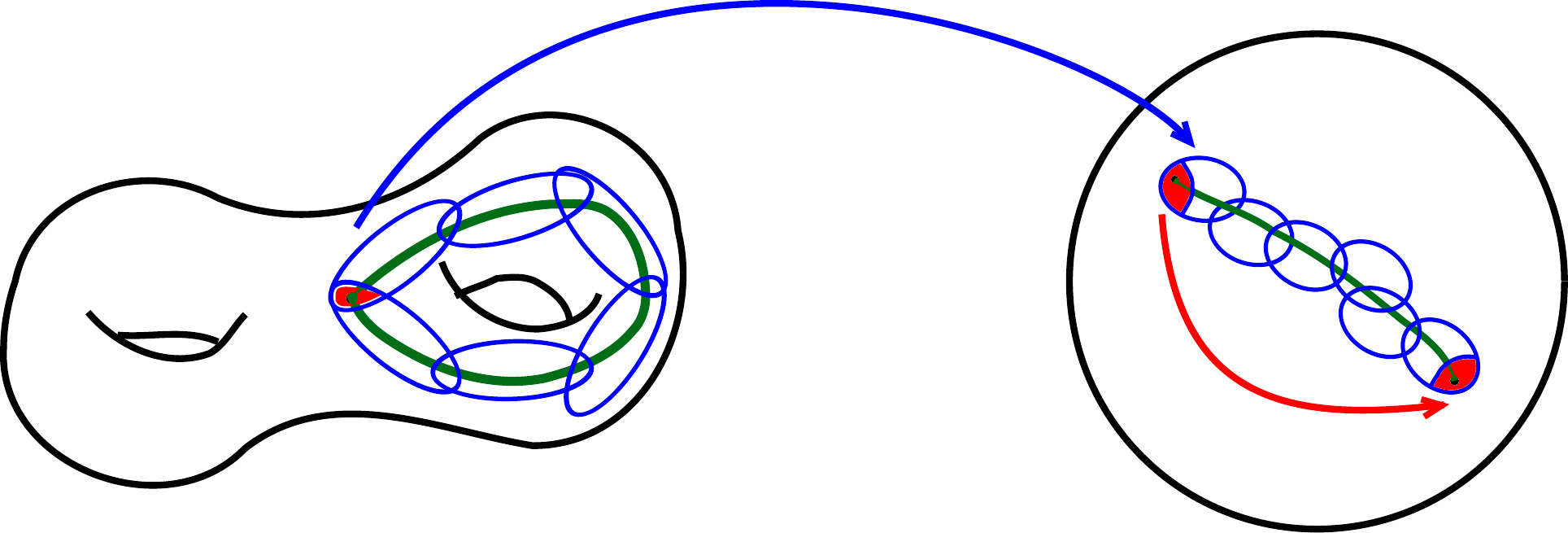}
\vspace{-0.5cm}
\caption{Construction of a holonomy representation}
\vspace{-0.4cm}
\label{fig-hol}
\end{figure}
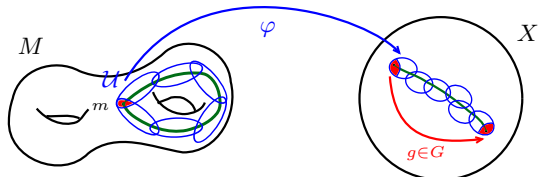

Let $M$ be a $(G,X)$-manifold, \ie a manifold endowed with a $(G,X)$-structure.
Fix a basepoint $m\in M$ and a chart $\varphi : \mathcal{U}\to X$ with $m\in\nolinebreak\mathcal{U}$.
We can lift any loop on~$M$ starting at~$m$ to a path on~$X$ starting at $\varphi(m)$, using successive charts of~$M$ which coincide on their intersections; the last chart in this analytic continuation process coincides, on an open set, with $g\cdot\varphi$ for some unique $g\in G$; we set $\mathrm{hol}(\gamma):=g$ where $\gamma\in\pi_1(M,m)$ is the homotopy class of the loop (see Figure~\ref{fig-hol}).
This defines a representation $\mathrm{hol} : \pi_1(M)\to G$ called the \emph{holonomy} (see \cite{gol-icm10,gol-book} for details); it is unique modulo conjugation by~$G$.
This coincides with the notion from Section~\ref{subsec:intro-GXstruct}; in particular, if $M\simeq\Gamma\backslash\mathcal{V}$ with $\mathcal{V}$ open in~$X$ and $\Gamma$ discrete in~$G$, and if $\mathcal{V}$ is simply connected, then $\mathrm{hol} : \pi_1(M)\to\Gamma$ is just the natural identification of $\pi_1(M)$ with~$\Gamma$.

We shall define the \emph{deformation space} $\mathrm{Def}_{(G,X)}(M)$ to be the quotient of the set of \emph{marked} $(G,X)$-structures on~$M$ (\ie pairs $(M',f')$ where $M'$ is a $(G,X)$-manifold and $f : M\to M'$ a diffeomorphism) by the group of diffeomorphisms of~$M$ isoto\-pic to the identity (acting by precomposition).
The holonomy defines a map from $\mathrm{Def}_{(G,X)}(M)$ to the space $\Hom(\Gamma,G)/G$ of representations of $\Gamma$ to~$G$~modulo conjugation by~$G$.
This map may be bijective in some cases, as in Example~\ref{ex:holonomy-Teich} below, but in general it is not.
However, when $M$ is closed, the so-called \emph{Eh\-resmann--Thurston principle} \cite{thu80} states that the map is continuous (for the natural topologies on both sides), open, with discrete fibers; in particular, the set of holonomy representations of $(G,X)$-structures on~$M$ is then stable under small~deformations.

\begin{example} \label{ex:holonomy-Teich}
Let $(G,X)=(\PO(2,1),\HH^2)$ where $\PO(2,1)\simeq\PGL(2,\R)$ is the isometry group of the real hyperbolic plane~$\HH^2$.
Let $M=S$ be a closed oriented connected surface of genus $g\geq 2$.
All $(G,X)$-structures on~$S$ are complete.
Their holonomy representations are the Fuchsian (\ie faithful and discrete) repre\-sentations from $\pi_1(S)$ to~$G$.
The deformation space $\mathrm{Def}_{(G,X)}(S)$ is the Teichm\"uller~space $\mathrm{Teich}(S)$.
The holonomy defines a homeomorphism between $\mathrm{Teich}(S) \simeq \R^{6g-6}$ and the space of Fuchsian representations from $\pi_1(S)$ to~$G$ modulo conjugation by~$G$.
\end{example}

\section{Examples of $(G,X)$-structures and their holonomy representations} \label{sec:dyn-properties-hol}

In this section we introduce three important families of $(G,X)$-structures, which have been much studied in the past few decades.
We observe some structural stability for their holonomy representations, and the existence of continuous equivariant boundary maps together with expansion/contraction properties ``at infinity''.
These phenomena will be captured by the notion of an Anosov representation in Section~\ref{sec:Anosov}.

\subsection{Convex cocompact locally symmetric structures in rank one} \label{subsec:cc-rank1}

Let $G$ be a real semisimple Lie group of real rank one with Riemannian symmetric space $X=G/K$ (\ie $K$ is a maximal compact subgroup of~$G$).
E.g.\ $(G,X)=(\PO(n,1),\HH^n)$ for $n\geq 2$.
Convex cocompact groups are an important class of discrete subgroups $\Gamma$ of~$G$ which generalize the uniform lattices.
They are special cases of geometrically finite groups, for which no cusps appear; see Bowditch \cite{bow93,bow98} for a general theory.

By definition, a discrete subgroup $\Gamma$ of~$G$ is \emph{convex cocompact} if it preserves and acts with compact quotient on some nonempty convex subset $\mathcal{C}$ of $X=G/K$; equivalently, the complete $(G,X)$-manifold (or orbifold) $\Gamma\backslash X$ has a compact convex subset (namely $\Gamma\backslash\mathcal{C}$) containing all the topology.
Such a group $\Gamma$ is always finitely~ge\-nerated.
A representation $\rho : \Gamma\to G$ is called \emph{convex cocompact} if its kernel~is~finite and its image is a convex cocompact subgroup of~$G$.

For instance, in Example~\ref{ex:qF-intro} the quasi-Fuchsian representations are exactly the convex cocompact representations from $\pi_1(S)$ to $G=\PSL(2,\C)$; modulo conjugation, they are parametrized by $\mathrm{Teich}(S)\times\mathrm{Teich}(S)$ \cite{ber60}.
Another classical example of convex cocompact groups in rank-one~$G$ is Schottky groups, namely free groups defined by the so-called \emph{ping pong} dynamics of their generators in $\partial_{\infty}X$.

Here $\partial_{\infty}X$ denotes the visual boundary of~$X$, yielding the standard compactification $\overline{X}=X\sqcup\partial_{\infty}X$ of~$X$; for $X=\HH^n$ we can see $\overline{X}$ in projective space as in Example~\ref{ex:conv-div}.(1) below.
The $G$-action on~$X$ extends continuously to~$\overline{X}$, and $\partial_{\infty}X$ identifies with $G/P$ where $P$ is a minimal parabolic subgroup of~$G$.

For a convex cocompact representation $\rho : \Gamma\to G$, the existence of a cocompact invariant convex set $\mathcal{C}$ implies (by the \v{S}varc--Milnor lemma or ``fundamental observation of geometric group theory'') that $\rho$ is a \emph{quasi-isometric embedding}.
This means that the points of any $\rho(\Gamma)$-orbit in $X=G/K$ go to infinity at linear speed for the word length function $|\cdot| : \Gamma\to\N$: for any $x_0\in X$ there exist $C,C'>0$ such that $d_{G/K}(x_0,\rho(\gamma)\cdot x_0) \geq C\,|\gamma| - C'$ for all $\gamma\in\Gamma$.
(This property does not depend on the choice of finite generating subset of $\Gamma$ defining~$|\cdot|$.)
A consequence ``at infinity'' is that any $\rho$-orbital map $\Gamma\to X$ extends to a $\rho$-equivariant embedding $\xi : \partial_{\infty}\Gamma\to\partial_{\infty}X\simeq G/P$, where $\partial_{\infty}\Gamma$ is the boundary of the Gromov hyperbolic group~$\Gamma$.
The image of~$\xi$ is the \emph{limit set} $\Lambda_{\rho(\Gamma)}$ of $\rho(\Gamma)$ in $\partial_{\infty}X$.
The dynamics on $\partial_{\infty}X\simeq G/P$ is decomposed as in Example~\ref{ex:qF-intro}: the action of $\rho(\Gamma)$ is ``chaotic'' on $\Lambda_{\rho(\Gamma)}$ (\eg all orbits are dense if $\Gamma$ is nonelementary), and properly discontinuous, with compact quotient, on the complement $\Omega_{\rho(\Gamma)}=\partial_{\infty}X\smallsetminus\Lambda_{\rho(\Gamma)}$.

Further dynamical properties were studied by Sullivan: by \cite{sul79}, the action of $\rho(\Gamma)$ on $\partial_{\infty}X\simeq G/P$ is \emph{expanding} at each point $z\in\Lambda_{\rho(\Gamma)}$, \ie there exist $\gamma\in\Gamma$ and $C>1$ such that $\rho(\gamma)$ multiplies all distances by $\geq C$ on a neighborhood of $z$ in $\partial_{\infty}X$ (for some fixed auxiliary metric on $\partial_{\infty}X$).
This implies that the group $\rho(\Gamma)$ is \emph{structurally stable}, \ie there is a neighborhood of the natural inclusion in $\Hom(\rho(\Gamma),G)$ consisting entirely of faithful representations.
In fact, $\rho$ admits a neighborhood consisting entirely of convex cocompact representations, by a variant of the Ehresmann--Thurston principle.
For $G=\SL(2,\C)$, a structurally stable subgroup of~$G$ is either locally rigid or convex cocompact, by \cite{sul85}.

\subsection{Convex projective structures: divisible convex sets} \label{subsec:conv-div}

Let $G$ be the projective linear group $\PGL(d,\R)$ and $X$ the projective space~$\PP(\R^d)$, for $d\geq 2$.
Recall that a subset of $X=\PP(\R^d)$ is said to be \emph{convex} if it is contained and convex in some affine chart, \emph{properly convex} if its closure is convex, and \emph{strictly convex} if it is properly convex and its boundary in~$X$ does not contain any nontrivial segment.

\begin{remark} \label{rem:Hilbert}
Any properly convex open subset $\Omega$ of $X=\PP(\R^d)$ admits a well-behaved (complete, proper, Finsler) metric $d_{\Omega}$, the \emph{Hilbert metric}, which is invariant under the subgroup of $G=\PGL(d,\R)$ preserving~$\Omega$ (see \eg \cite{ben-cd-survey}).
In particular, any discrete subgroup of~$G$ preserving~$\Omega$ acts properly discontinuously on~$\Omega$.
\end{remark}

By definition, a \emph{convex projective structure} on a manifold~$M$ is a $(G,X)$-structure obtained by identifying $M$ with $\Gamma\backslash\Omega$ for some properly convex open subset $\Omega$ of~$X$ and some discrete subgroup $\Gamma$ of~$G$.
When $M$ is closed, \ie when $\Gamma$ acts with compact quotient, we say that $\Gamma$ \emph{divides}~$\Omega$.
Such \emph{divisible convex sets} $\Omega$ are the objects of a rich theory, see \cite{ben-cd-survey}.
The following classical examples are called \emph{symmetric}.

\begin{examples} \label{ex:conv-div}
(1) For $d=n+1\geq 3$, let $\langle\cdot,\cdot\rangle_{n,1}$ be a symmetric bilinear form of signature $(n,1)$ on~$\R^d$, and $\Omega = \{ [v]\in\PP(\R^d) \,|\, \langle v,v\rangle_{n,1} < 0\}$ be the projective~mo\-del of the real hyperbolic space~$\HH^n$.
It is a strictly convex open subset of $X=\PP(\R^d)$ (an ellipsoid), and any uniform lattice $\Gamma$ of $\PO(n,1)\subset G=\PGL(d,\R)$ divides~$\Omega$.

(2) For $d=n(n+1)/2$, let us see $\R^d$ as the space $\mathrm{Sym}(n,\R)$ of symmetric $n\times n$ real matrices, and let $\Omega\subset\PP(\R^d)$ be the image of the set of positive definite ones.
The set $\Omega$ is a properly convex open subset of $X=\PP(\R^d)$; it is strictly convex if and only if $n=2$.
The group $\GL(n,\R)$ acts on $\mathrm{Sym}(n,\R)$ by $g\cdot s:=gs^t\!g$, which induces an action of $\PGL(n,\R)$~on~$\Omega$.~This~action is transitive and the stabilizer of a point is $\PO(n)$, hence $\Omega$ identifies with the Riemannian symmetric space $\PGL(n,\R)/\PO(n)$.
In particular, any uniform lattice $\Gamma$ of $\PGL(n,\R)$ divides~$\Omega$.
(A similar construction works over the complex numbers, the quaternions, or the octonions: see \cite{ben-cd-survey}.)
\end{examples}

Many nonsymmetric strictly examples were also constru\-cted since the 1960s by various techniques; see \cite{ben-cd-survey,clm16-survey} for references.
Remarkably, there exist irreducible di\-visible convex sets $\Omega\subset\PP(\R^d)$ which are not symmetric and not strictly convex: the first examples were built by Benoist \cite{ben-cd-IV} for $4\leq d\leq 7$.
Ballas--Danciger--Lee~\cite{bdl} generalized Benoist's construction for $d=4$ to show that large families of nonhyperbolic closed $3$-manifolds admit convex projective structures.
Choi--Lee--Mar\-quis \cite{clm16} recently built nonstrictly convex examples of a different~flavor~for~$5\leq\nolinebreak d\leq 7$.

For \emph{strictly convex}~$\Omega$, dynamics ``at infinity'' are relatively well understood: if $\Gamma$ divides~$\Omega$, then $\Gamma$ is Gromov hyperbolic \cite{ben-cd-I} and, by cocompactness, any orbital map $\Gamma\to\Omega$ extends continuously to an equivariant homeomorphism from the boundary $\partial_{\infty}\Gamma$ of~$\Gamma$ to the boundary of $\Omega$ in~$X$.
This is similar to Section~\ref{subsec:cc-rank1}, except that now $X$ itself is a flag variety $G/P$ (see Table~\ref{table1}).
The image of the boundary map is again a \emph{limit set} $\Lambda_{\Gamma}$ on which the action of~$\Gamma$ is ``chaotic'', but $\Lambda_{\Gamma}$ is now part of a larger ``extended limit set'' $\mathscr{L}_{\Gamma}$, namely the union of all projective hyperplanes tangent to~$\Omega$ at points of $\Lambda_{\Gamma}$.
The space $X\simeq G/P$ is the disjoint union of $\mathscr{L}_{\Gamma}$ and~$\Omega$.
The dynamics of $\Gamma$ on $X$ are further understood by considering the \emph{geodesic flow} on $\Omega\subset X$, defined using the Hilbert metric of Remark~\ref{rem:Hilbert}; for $\Omega=\HH^n$ as in Example~\ref{ex:conv-div}.(1), this is the usual geodesic flow.
Benoist \cite{ben-cd-I} proved that the induced flow on $\Gamma\backslash\Omega$ is Anosov and topologically mixing; see \cite{cra14} for further~properties.

For \emph{nonstrictly convex}~$\Omega$, the situation is less understood.
Groups $\Gamma$ dividing~$\Omega$ are never Gromov hyperbolic \cite{ben-cd-I}; for $d=4$ they are relatively hyperbolic \cite{ben-cd-IV}, but in general they might not be (\eg if $\Omega$ is symmetric), and it is not obvious what~type of boundary $\partial_{\infty}\Gamma$ (defined independently of~$\Omega$) might be most~useful~in~the context of Problem~\ref{problemB}.
The geodesic flow on $\Gamma\backslash\Omega$ is not Anosov, but Bray~\cite{bra}~proved it is still topologically mixing for $d=4$.
Much of the dynamics~remains~to~be~explored.

By Koszul \cite{kos68}, discrete subgroups of~$G$ dividing~$\Omega$ are structurally stable; moreover, for a closed manifold~$M$ with fundamental group $\Gamma=\pi_1(M)$, the set $\Hom_M^{\mathrm{conv}}(\Gamma,G)$ of holonomy representations of convex $(G,X)$-structures on~$M$ is open in $\Hom(\Gamma,G)$.
It is also closed in $\Hom(\Gamma,G)$ as soon as $\Gamma$ does not contain an infinite normal abelian subgroup, by Choi--Goldman \cite{cg05} (for $d=3$)~and~Benoist~\cite{ben-cd-III} (in general).
For $d=3$, when $M$ is a closed surface of genus $g\geq 2$, Goldman \cite{gol90}\linebreak showed that $\Hom_M^{\mathrm{conv}}(\Gamma,G)/G$ is homeomorphic to $\R^{16g-16}$, via an explicit~parame\-trization generalizing classical (\emph{Fenchel--Nielsen}) coordinates on Teichm\"uller space.%

\subsection{AdS quasi-Fuchsian representations} \label{subsec:AdS-qF}

We now discuss the Lorentzian counterparts of Example~\ref{ex:qF-intro}, which have been studied by Witten \cite{wit88} and others as simple models for $(2+1)$-dimensional gravity.
Let $M=S\times (0,1)$ be as in Exam\-ple~\ref{ex:qF-intro}.
Instead of taking $(G,X)=(\PO(3,1),\HH^3)$, we now take $G=\PO(2,2)$~and
$$X = \AdS^3 = \{ [v]\in\PP(\R^4) ~|~ \langle v,v\rangle_{2,2}<0\}.$$
In other words, we change the signature of the quadratic form defining~$X$ from $(3,1)$ (as in Example~\ref{ex:conv-div}.(1)) to $(2,2)$.
This changes the natural $G$-invariant metric from Riemannian to Lorentzian, and the topology of~$X$ from a ball to a solid torus.
The space $X=\AdS^3$ is called the \emph{anti-de Sitter} $3$-space.

The manifold $M=S\times (0,1)$ does not admit $(G,X)$-structures of type~C (see Section~\ref{subsec:intro-GXstruct}), but it admits some of type~U, called \emph{globally hyperbolic~maxi\-mal Cauchy-compact} (GHMC).
In general, a Lorentzian manifold is called globally hyperbolic if it satisfies the intuitive property that ``when moving towards the future one does not come back to the past''; more precisely, there is a spacelike hypersurface (\emph{Cauchy hypersurface}) meeting each inextendible causal curve exactly once.
Here we also require that the Cauchy surface be compact and that $M$ be maximal (\ie not isometrically embeddable into a larger globally hyperbolic Lorentzian $3$-manifold).

To describe the GHMC $(G,X)$-structures on~$M$, it is convenient to consider a different model for $\AdS^3$, which leads to beautiful links with $2$-dimensional hyperbolic geometry.
Namely, we view $\R^4$ as the space $\mathcal{M}_2(\R)$ of real $2\times 2$ matrices, and the quadratic form $\langle\cdot,\cdot\rangle_{2,2}$ as minus the determinant.
This induces an identification of $X=\AdS^3$ with $\underline{G}=\PSL(2,\R)$ sending $[v]\in X$ to $\big[ \frac{1}{|\langle v,v\rangle|}(\begin{smallmatrix} v_1+v_4 & v_2+v_3\\ v_2-v_3 & -v_1+v_4\end{smallmatrix}) \big] \in \underline{G}$, and a corresponding group isomorphism from the identity component $G_0=\PO(2,2)_0$ of~$G$ acting on $X=\AdS^3$, to $\underline{G}\times\underline{G}$ acting on $\underline{G}$ by right and left multiplication:~$(g_1,g_2)\cdot g=g_2gg_1^{-1}$.
It also induces an identification of the boundary $\partial X\subset\PP(\R^4)$ with the projectivization of the set of rank-one matrices, hence with $\PP^1(\R)\times\PP^1(\R)$ (by taking the kernel and the image); the action of $G_0$ on $\partial X$ corresponds to the natural action of $\underline{G}\times\underline{G}$~on~$\PP^1(\R)\times\PP^1(\R)$.

\begin{figure}[ht!]
\labellist
\small\hair 2pt
\pinlabel {$\AdS^3$} [r] at 40 105
\pinlabel {$\Omega$} [r] at 70 90
\pinlabel {$\Lambda$} [r] at 78 62
\pinlabel {$\mathcal{C}$} [r] at 78 53
\endlabellist
\centering
\vspace{-0.3cm}
\includegraphics[scale=1]{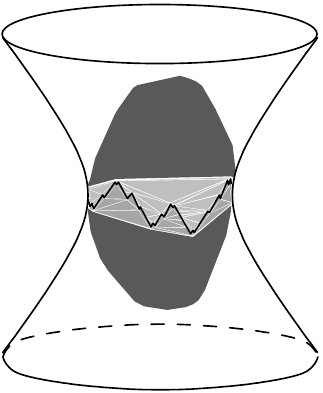}
\vspace{-0.2cm}
\caption{The sets $\Lambda$, $\Omega$, $\mathcal{C}$ for an AdS quasi-Fuchsian representation}
\vspace{-0.2cm}
\label{fig-AdS3}
\end{figure}
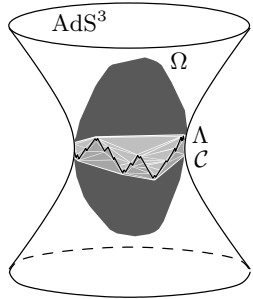

With these identifications, Mess \cite{mes90} proved that all GHMC $(G,X)$-structures~on $M=S\times\nolinebreak (0,1)$ are obtained as follows.
Let $(\rho_L,\rho_R)$ be a pair of Fuchsian representations from $\Gamma=\pi_1(M)\simeq\pi_1(S)$ to $\underline{G}=\PSL(2,\R)$.
The group $(\rho_L,\rho_R)(\Gamma) \subset \underline{G}\times\underline{G} \subset G$ preserves a topological circle $\Lambda$ in $\partial X$, namely the graph of the homeomor\-phism of $\PP^1(\R)$ conjugating the action of $\rho_L$ to that of~$\rho_R$.
For any $z\in\Lambda$, the orthogonal $z^{\perp}$ of~$z$ for $\langle\cdot,\cdot\rangle_{2,2}$ is a projective hyperplane tangent to $X$ at~$z$.
The complement $\Omega$ in $\PP(\R^4)$ of the union of all $z^{\perp}$ for $z\in\Lambda$ is a convex open subset of $\PP(\R^4)$ contained in~$X$ (see Figure~\ref{fig-AdS3}) which admits a $\Gamma$-invariant Cauchy surface.
The action of $\Gamma$ on $\Omega$ via $(\rho_L,\rho_R)$ is properly discontinuous and the convex hull $\mathcal{C}$ of $\Lambda$ in~$\Omega$ (called the \emph{convex core}) has compact quotient by~$\Gamma$.
The quotient $(\rho_L,\rho_R)(\Gamma)\backslash\Omega$ is diffeomorphic to $M=S\times (0,1)$, and this yields a GHMC $(G,X)$-structure on~$M$.

Such $(G,X)$-structures, or their holonomy representations $\rho = (\rho_L,\rho_R) : \Gamma\to\underline{G}\times\underline{G}\subset G$, are often called \emph{AdS quasi-Fuchsian}, by analogy with Example~\ref{ex:qF-intro}.
Their deformation space
is parame\-trized by $\mathrm{Teich}(S)\times\mathrm{Teich}(S)$, via $(\rho_L,\rho_R)$ \cite{mes90}.
Their geometry, especially the geometry of the convex core and the way it determines $(\rho_L,\rho_R)$, is the object of active current research (see \cite{bbd+,bs12}).
Generalizations have recently been worked out in several directions (see \cite{bar15,bm12,bks11} and Section~\ref{subsec:cc-Hpq}).

As in Section~\ref{subsec:cc-rank1}, the compactness of the convex core of an AdS quasi-Fuchsian manifold implies that any orbital map $\Gamma\to\Omega$ extends ``at infinity'' to an equivariant embedding $\xi : \partial_{\infty}\Gamma\to\partial X$ with image~$\Lambda$.
Here $\partial X$ is still a flag variety $G/P$, where $P$ is the stabilizer in $G=\PO(2,2)$ of an isotropic line of $\R^4$ for $\langle\cdot,\cdot\rangle_{2,2}$.
Although $G$ has higher rank, the rank-one dynamics of Section~\ref{subsec:cc-rank1} still appear through the product structure of $G_0 \simeq \underline{G}\times\underline{G}$ acting on $\partial X \simeq \PP^1(\R)\times\PP^1(\R) \simeq \partial_{\infty}\HH^2\times\partial_{\infty}\HH^2$.

\section{Anosov representations} \label{sec:Anosov}

In this section we define and discuss Anosov representations.
These are representations of Gromov hyperbolic groups into Lie groups~$G$ with strong dynamical properties, defined using continuous equivariant boundary maps.
They were introduced by Labourie \cite{lab06} and further investigated by~Gui\-chard--Wienhard \cite{gw12}.
They play an important role in higher~Teich\-m\"uller theory and in the study of Problem~\ref{problemB}.
As we shall see in Section~\ref{subsec:Section-3-was-Ano}, most representations that appeared in Section~\ref{sec:dyn-properties-hol} were in fact Anosov representations.

\subsection{The definition} \label{subsec:def-Ano}

Let $\Gamma$ be a Gromov hyperbolic group with boundary $\partial_{\infty}\Gamma$ (\eg $\Gamma$ a surface group and $\partial_{\infty}\Gamma$ a circle, or $\Gamma$ a nonabelian free group and $\partial_{\infty}\Gamma$ a Cantor set).
The notion of an Anosov representation of $\Gamma$ to a reductive Lie group~$G$ depends on the choice of a parabolic subgroup $P$ of~$G$ up to conjugacy, \ie on the choice of a flag variety $G/P$ (see Section~\ref{subsec:intro-G/P}).
Here, for simplicity, we restrict to $G=\PGL(d,\R)$.
We choose an integer $i\in [1,d-1]$ and denote by $P_i$ the stabilizer in~$G$ of an $i$-plane of~$\R^d$, so that $G/P_i$ identifies with the Grassmannian~$\mathrm{Gr}_i(\R^d)$.

By definition, a representation $\rho : \Gamma\to\PGL(d,\R)$ is \emph{$P_i$-Anosov} if there exist two con\-tinuous $\rho$-equivariant maps $\xi_i : \partial_{\infty}\Gamma\to\mathrm{Gr}_i(\R^d)$ and $\xi_{d-i} :\nolinebreak\partial_{\infty}\Gamma\to\nolinebreak\mathrm{Gr}_{d-i}(\R^d)$ which are transverse (\ie $\xi_i(\eta)+\xi_{d-i}(\eta')=\R^d$ for all $\eta\neq\eta'$ in $\partial_{\infty}\Gamma$) and satisfy a uniform contraction/expansion condition analogous to that defining Anosov flows.

Let us state this condition in the original case considered by Labourie \cite{lab06}, where $\Gamma=\pi_1(M)$ for some closed negatively-curved manifold~$M$.
We denote by $\widetilde{M}$ the universal cover of~$M$, by $T^1$ the unit tangent bundle, and by $(\varphi_t)_{t\in\R}$ the geodesic flow on either $T^1(M)$ or $T^1(\widetilde{M})$.~Let
$$E^{\rho} = \Gamma\backslash (T^1(\widetilde{M})\times\R^d)$$
be the natural flat vector bundle over $T^1(M)=\Gamma\backslash T^1(\widetilde{M})$ associated to~$\rho$, where $\Gamma$ acts on $T^1(\widetilde{M})\times\R^d$ by $\gamma\cdot (\tilde{x},v) = (\gamma\cdot\tilde{x}, \rho(\gamma)\cdot v)$.
The geodesic flow $(\varphi_t)_{t\in\R}$ on $T^1(M)$ lifts to a flow $(\psi_t)_{t\in\R}$ on $E^{\rho}$, given by $\psi_t\cdot [(\tilde{x},v)] = [(\varphi_t\cdot\tilde{x},v)]$.
For each $\tilde{x}\in\nolinebreak T^1(\widetilde{M})$, the transversality of the boundary maps induces a decomposition $\R^d = \xi_i(\tilde{x}^+) \oplus \xi_{d-i}(\tilde{x}^-)$, where $\tilde{x}^{\pm} = \lim_{t\to\pm\infty} \varphi_t\cdot\tilde{x}$ are the forward and backward endpoints of the geodesic defined by~$\tilde{x}$, and this defines a decomposition of the vector bundle $E^{\rho}$ into the direct sum of two subbundles $E_i^{\rho} = \{ [(\tilde{x},v)] \,|\, v\in\xi_i(\tilde{x}^+)\}$ and $E_{d-i}^{\rho} = \{ [(\tilde{x},v)] \,|\, v\in\xi_{d-i}(\tilde{x}^-)\}$.
This decomposition is invariant under the flow $(\psi_t)$.
By definition, the representation $\rho$ is $P_i$-Anosov if the following condition~is~satisfied.

\begin{condition} \label{cond:Ano}
The flow $(\psi_t)_{t\in\R}$ uniformly contracts $E_i^{\rho}$ with respect to $E_{d-i}^{\rho}$, \ie there exist $C,C'>0$ such that for any $t\geq 0$, any $x\in T^1(M)$, and any nonzero $w_i\in E_i^{\rho}(x)$ and $w_{d-i}\in E_{d-i}^{\rho}(x)$,
$$\frac{\Vert\psi_t\cdot w_i\Vert_{\varphi_t\cdot x}}{\Vert\psi_t\cdot w_{d-i}\Vert_{\varphi_t\cdot x}} \leq e^{-Ct+C'} \, \frac{\Vert w_i\Vert_x}{\Vert w_{d-i}\Vert_x},$$
where $(\Vert\cdot\Vert_x)_{x\in T^1(M)}$ is any fixed continuous family of norms on the fibers $E^{\rho}(x)$.
\end{condition}

See \cite{bcls15} for an interpretation in terms of \emph{metric Anosov flows} (or Smale flows).

Condition~\ref{cond:Ano} implies in particular that the boundary maps $\xi_i$, $\xi_{d-i}$ are \emph{dynamics-preserving}, in the sense that the image of the attracting fixed point in $\partial_{\infty}\Gamma$ of any infinite-order element $\gamma\in\Gamma$ is an attracting fixed point in $\mathrm{Gr}_i(\R^d)$ or $\mathrm{Gr}_{d-i}(\R^d)$ of $\rho(\gamma)$.
Thus $\xi_i$ and~$\xi_{d-i}$ are unique, by density of such fixed points in $\partial_{\infty}\Gamma$.

We note that \emph{$P_i$-Anosov is equivalent to $P_{d-i}$-Anosov}, as the integers $i$ and $d-i$ play a similar role in the definition (up to reversing the flow, which switches contraction and expansion).
In particular, we may restrict to $P_i$-Anosov for $1\leq\nolinebreak i\leq\nolinebreak d/2$.

Guichard--Wienhard \cite{gw12} observed that an analogue of Condition~\ref{cond:Ano} can actually be defined for any Gromov hyperbolic group~$\Gamma$.
The idea is to replace $T^1(\widetilde{M})$ by $\partial_{\infty}\Gamma^{(2)}\times\R$ where $\partial_{\infty}\Gamma^{(2)}$ is the space of pairs of distinct points in the Gromov boundary $\partial_{\infty}\Gamma$ of~$\Gamma$, and the flow $\varphi_t$ by translation by $t$ along the $\R$ factor.
The work of Gromov \cite{gro87} (see also \cite{cha94,mat91,min05}) yields an appropriate extension of the $\Gamma$-action on $\partial_{\infty}\Gamma^{(2)}$ to $\partial_{\infty}\Gamma^{(2)}\times\R$, which is properly discontinuous and cocompact.
This leads to a notion of an Anosov representation for any Gromov hyperbolic~group~$\Gamma$~\cite{gw12}.

\subsection{Important properties and examples} \label{subsec:ex-Anosov}

A fundamental observation motivating the study of Anosov representations is the following: if $G$ is a semisimple Lie group of real rank one, then a representation $\rho : \Gamma\to G$ is Anosov if and only if it is convex cocompact in the sense of Section~\ref{subsec:cc-rank1}.

Moreover, many important properties of convex cocompact representations into rank-one groups generalize to Anosov representations.
For instance, Anosov representations $\rho : \Gamma\to G$ are quasi-isometric embeddings \cite{gw12,lab06}; in particular, they have finite kernel and discrete image.
Also by \cite{gw12,lab06}, any Anosov subgroup (\ie the image of any Anosov representation $\rho : \Gamma\to G$) is structurally stable; moreover, $\rho$ admits a neighborhood in $\Hom(\Gamma,G)$ consisting entirely of Anosov representations.
This is due to the uniform hyperbolicity nature of the Anosov condition.

Kapovich, Leeb, and Porti, in a series of papers (see~\cite{kl17,klp-survey,gui-bourbaki}), have developed a detailed analogy between Anosov representations into higher-rank semisimple Lie groups and convex cocompact representations into rank-one simple groups, from the point of view of dynamics (\eg extending the expansion property at the~limit~set~of Section~\ref{subsec:cc-rank1} and other classical characterizations) and\,topology\,(e.g.\,compactifications).

Here are some classical examples of Anosov representations in higher real rank.

\begin{examples} \label{ex:classical-Ano}
Let $\Gamma=\pi_1(S)$ where $S$ is a closed orientable surface of genus $\geq 2$.

(1) (Labourie \cite{lab06}) For $d\geq 2$, let $\tau_d : \PSL(2,\R)\to G=\PGL(d,\R)$ be the irreducible representation (unique up to conjugation by~$G$).
For any Fuchsian representation $\rho_0 :\nolinebreak\Gamma\to\PSL(2,\R)$, the composition $\tau_d\circ\rho_0 : \Gamma\to G$ is $P_i$-Anosov for all $1\leq i\leq d-1$.
Moreover, any representation in the connected component of $\tau_d\circ\rho_0$ in $\Hom(\Gamma,G)$ is still $P_i$-Anosov for all $1\leq i\leq d-1$.
These representations were first studied by Hitchin \cite{hit92} and are now known as \emph{Hitchin representations}.

(2) (Burger--Iozzi--Labourie--Wienhard \cite{bilw05,biw})
If a representation of $\Gamma$ into $G=\mathrm{PSp}(2n,\R)\subset\PGL(2n,\R)$ (\resp $G=\PO(2,q)\subset\PGL(2+q,\R)$) is \emph{maximal}, then it is $P_n$-Anosov (\resp $P_1$-Anosov).

(3) (Barbot \cite{bar10} for $d=3$) Let $d\geq 2$.
Any Fuchsian representation $\Gamma\to\SL(2,\R)$, composed with the standard embedding $\SL(2,\R)\hookrightarrow\SL(d,\R)$ (given by the direct sum of the standard action on~$\R^2$ and the trivial action on $\R^{d-2}$), defines a $P_1$-Anosov representation $\Gamma\to G=\PSL(d,\R)$.
\end{examples}

In~(2), we say that $\rho : \Gamma\to G$ is \emph{maximal} if it maximizes a topological invariant, the \emph{Toledo number} $T(\rho)$, defined for any simple $G$ of Hermitian type.
If $X=G/K$ is the Riemannian sym\-metric space of~$G$, then the imaginary part of the $G$-invariant Hermitian form on~$X$ defines a real $2$-form $\omega_X$, and by definition $T(\rho) = \frac{1}{2\pi} \, \int_S f^*\omega_X$ where $f : \widetilde{S}\to X$ is any~$\rho$-equi\-variant smooth map.
For $G=\PSL(2,\R)$, this coincides with the Euler number of~$\rho$.~In~ge\-neral, $T(\rho)$ takes discrete values and $|T(\rho)|\leq\mathrm{rank}_{\R}(G)\,|\chi(S)|$ where $\chi(S)$ is the Euler characteristic of~$S$ (see \cite{biw14}).

While (1) and~(3) provide Anosov representations in two of the three connected components of $\Hom(\Gamma,\PSL(3,\R))$ for $\Gamma=\pi_1(S)$, it is currently not known whether Anosov representations appear in the third component.

See \cite{ben96,bt17,cls17,gw12,klp14bis} for higher-rank Anosov generalizations of Schottky groups.

\subsection{Higher Teichm\"uller spaces of Anosov representations} \label{subsec:Higher-Teich}

Anosov representations play an important role in \emph{higher Teich\-m\"uller theory}, a currently very active theory whose goal is to find deformation spaces of faithful and discrete representations of discrete groups $\Gamma$ into higher-rank semisimple Lie groups $G$ which share some of the remarkable properties of Teichm\"uller space.
Although various groups $\Gamma$ may be considered, the foundational case is when $\Gamma=\pi_1(S)$ for some closed connected surface~$S$ of genus $\geq 2$ (see \cite{biw14,wie-icm}); then one can use rich features of Riemann surfaces, explicit topological considerations, and powerful techniques based on Higgs bundles as in Hitchin's pioneering~work~\cite{hit92}.

Strikingly similar properties to $\mathrm{Teich}(S)$ have been found for two types of \emph{higher Teichm\"uller spaces}: the space of \emph{Hitchin representations} of $\Gamma$ into a real split simple Lie group $G$ such as $\PGL(d,\R)$, modulo conjugation by~$G$; and the space of \emph{maximal representations} of $\Gamma$ into a simple Lie group $G$ of Hermitian type such as $\mathrm{PSp}(2n,\R)$ or $\PO(2,q)$, modulo conjugation by~$G$.
Both these spaces are unions of connected components of $\Hom(\Gamma,G)/G$, consisting entirely of Anosov representations (see Examples~\ref{ex:classical-Ano}.(1)--(2)).
Similarities of these spaces to $\mathrm{Teich}(S)$~include:
\begin{enumerate}
  \item the proper discontinuity of the action of the mapping class group $\mathrm{Mod}(S)$~\cite{lab08};
  \item for Hitchin representations to $\PGL(d,\R)$: the topology of $\R^{(d^2-1)\,|\chi(S)|}$~\cite{hit92};
  \item good\,systems\,of\,coordinates\,generalizing\,those\,on\,$\mathrm{Teich}(S)$\,\cite{bd14,fg06,gol90,str15,zha15};
  \item an\,analytic\,$\mathrm{Mod}(S)$-invariant\,Riemannian\,metric\,(\emph{pressure\,metric})\,\cite{bcls15,ps16};
  \item versions of the collar lemma for associated locally symmetric spaces~\cite{bp17,lz17}.
\end{enumerate}
Other higher Teichm\"uller spaces of Anosov representations of $\pi_1(S)$ are also being explored \cite{gw16}.
We refer to Section~\ref{sec:geom-for-Ano} for geometric structures associated~to~such~spaces.

\subsection{Characterizations} \label{subsec:charact-Ano}

Various characterizations of Anosov representations have been developed in the past few years, by Labourie \cite{lab06}, Guichard--Wienhard \cite{gw12}, Kapovich--Leeb--Porti \cite{kl17,klp14,klp14bis}, Gu\'eritaud--Guichard--Kassel--Wienhard \cite{ggkw17}, and others.
Here are some characterizations that do not involve any flow.
They hold for any reductive Lie group~$G$, but for simplicity we state them for $G=\PGL(d,\R)$.
For $1\leq i\leq d$ and $g\in\GL(d,\R)$, we denote by $\mu_i(g)$ (\resp $\lambda_i(g)$) the logarithm of the $i$-th singular value (\resp eigenvalue) of~$g$.

\begin{theorem} \label{thm:charact-Ano}
For a Gromov hyperbolic group~$\Gamma$, a representation $\rho : \Gamma\to G=\PGL(d,\R)$, and an integer $1\leq i\leq d/2$, the following are equivalent:
\begin{enumerate}[leftmargin=0.8cm]
  \item\label{item:Ano} $\rho$ is $P_i$-Anosov {\small (or equivalently $P_{d-i}$-Anosov, see Section~\ref{subsec:def-Ano})};
  \item\label{item:Ano-xi-mu} there exist continuous, $\rho$-equivariant, transverse, dynamics-preserving boundary maps $\xi_i : \partial_{\infty}\Gamma\to\mathrm{Gr}_i(\R^d)$ and $\xi_{d-i} : \partial_{\infty}\Gamma\to\mathrm{Gr}_{d-i}(\R^d)$, and\\ $(\mu_i-\mu_{i+1})(\rho(\gamma))\to +\infty$ as $|\gamma|\to +\infty$;
  \item\label{item:Ano-xi-lambda} there exist continuous, $\rho$-equivariant, transverse, dynamics-preserving boundary maps $\xi_i : \partial_{\infty}\Gamma\to\mathrm{Gr}_i(\R^d)$ and $\xi_{d-i} : \partial_{\infty}\Gamma\to\mathrm{Gr}_{d-i}(\R^d)$, and\\ $(\lambda_i-\lambda_{i+1})(\rho(\gamma))\to +\infty$ as $\ell_{\Gamma}(\gamma)\to +\infty$;
  \item\label{item:Ano-mu} there exist $C,C'>0$ such that $(\mu_i-\mu_{i+1})(\rho(\gamma))\geq C\,|\gamma|-C'$ for all $\gamma\in\Gamma$;
  \item\label{item:Ano-lambda} there exist $C,C'>0$ such that $(\lambda_i-\lambda_{i+1})(\rho(\gamma))\geq C\,\ell_{\Gamma}(\gamma)-C'$ for all $\gamma\in\Gamma$.
\end{enumerate}
\end{theorem}

Here we denote by $|\cdot| : \Gamma\to\N$ the word length with respect to some fixed finite generating subset  of~$\Gamma$, and by $\ell_{\Gamma} : \Gamma\to\N$ the translation length in the Cayley graph of $\Gamma$ for that subset, \ie $\ell_{\Gamma}(\gamma) = \min_{\beta\in\Gamma} |\beta\gamma\beta^{-1}|$.
In a Gromov hyperbolic group~$\Gamma$ the translation length $\ell_{\Gamma}(\gamma)$ is known to differ only by at most a uniform additive constant from the stable length $|\gamma|_{\infty} = \lim_{n\to +\infty} |\gamma^n|/n$, and so we may replace $\ell_{\Gamma}(\gamma)$ by $|\gamma|_{\infty}$ in Conditions \eqref{item:Ano-xi-lambda} and~\eqref{item:Ano-lambda}.

The equivalence \eqref{item:Ano} $\Leftrightarrow$ \eqref{item:Ano-xi-mu} is proved in \cite{ggkw17} and \cite{klp14}, the equivalence \eqref{item:Ano-xi-mu} $\Leftrightarrow$ \eqref{item:Ano-xi-lambda} in \cite{ggkw17}, the equivalence \eqref{item:Ano} $\Leftrightarrow$ \eqref{item:Ano-mu} in \cite{klp14bis} and \cite{bps}, and the equivalence \eqref{item:Ano-mu} $\Leftrightarrow$ \eqref{item:Ano-lambda}~in~\cite{kp17}.

Condition~\eqref{item:Ano-mu} is a refinement of the condition of being a quasi-isometric embedding, which for $G=\PGL(d,\R)$ is equivalent to the existence of $C,C'>0$ such that $\sqrt{\sum_k (\mu_k-\mu_{k+1})^2(\rho(\gamma))} \geq C\,|\gamma| - C'$ for all $\gamma\in\Gamma$.
We refer to \cite{ggkw17} (\emph{CLI~con\-dition}) or \cite{klp14bis} (\emph{Morse condition}) for further refinements\,satisfied\,by\,Anosov\,representations.

By \cite{klp14bis,bps}, if $\Gamma$ is any finitely generated group, then the existence of a representation $\rho : \Gamma\to\PGL(d,\R)$ satisfying Condition~\eqref{item:Ano-mu} implies that $\Gamma$ is Gromov hyperbolic.
The analogue for~\eqref{item:Ano-lambda} is more subtle: \eg the Baumslag--Solitar group $\mathrm{BS}(1,2)$, which is not Gromov hyperbolic, still admits a faithful representation into $\PSL(2,\R)$ satisfying Condition~\eqref{item:Ano-lambda} for the stable length $|\cdot|_{\infty}$, see~\cite{kp17}.

Kapovich--Leeb--Porti's original proof \cite{klp14bis} of \eqref{item:Ano}~$\Leftrightarrow$~\eqref{item:Ano-mu} uses the geometry of higher-rank Riemannian symmetric spaces and asymptotic cones.
Bochi--Potrie--Sambarino's alternative proof \cite{bps} is based on an interpretation of \eqref{item:Ano} and~\eqref{item:Ano-mu} in terms of partially hyperbolic dynamics, and more specifically of dominated splittings for locally constant linear cocycles over certain subshifts.
Pursuing this point of view further, \cite{kp17} shows that the equivalence~\eqref{item:Ano-mu}~$\Leftrightarrow$~\eqref{item:Ano-lambda} of Theorem~\ref{thm:charact-Ano} implies the equivalence between nonuniform hyperbolicity (\ie all invariant measures are hyperbolic) and uniform hyperbolicity for a certain cocycle naturally associated with~$\rho$ on the space of biinfinite geodesics of~$\Gamma$.
In general in smooth dynamics, nonuniform hyperbolicity does not imply uniform hyperbolicity.

\subsection{Revisiting the examples of Section~\ref{sec:dyn-properties-hol}} \label{subsec:Section-3-was-Ano}

The boundary maps and dynamics ``at infinity'' that appeared in most examples of Section~\ref{sec:dyn-properties-hol} are in fact explained by the notion of an Anosov representation:
\begin{itemize}[leftmargin=0.8cm]
  \item convex cocompact representations into rank-one simple Lie groups as in Section~\ref{subsec:cc-rank1} are all Anosov (see Section~\ref{subsec:ex-Anosov});
  \item if $S$ is a closed orientable connected surface of genus $\geq 2$, then by Choi--Goldman \cite{cg05,gol90} the holonomy representations of convex projective structures on~$S$ as in Section~\ref{subsec:conv-div} are exactly the Hitchin representations of $\pi_1(S)$ into $\PSL(3,\R)$; they are all $P_1$-Anosov (Example~\ref{ex:classical-Ano}.(1));
  \item for general $d\geq 3$, Benoist's work \cite{ben-cd-I} shows that if $\Gamma$ is a discrete subgroup of $\PGL(d,\R)$ dividing a \emph{strictly convex} open subset $\Omega$ of $\PP(\R^d)$, then $\Gamma$ is Gromov hyperbolic and the inclusion $\Gamma\hookrightarrow\PGL(d,\R)$~is~$P_1$-Anosov;
  \item Mess's theory \cite{mes90} implies that a representation $\rho : \pi_1(S)\to\PO(2,2)\subset\PGL(4,\R)$ is AdS quasi-Fuchsian if and only if it is $P_1$-Anosov.
\end{itemize}

\section{Geometric structures for Anosov representations} \label{sec:geom-for-Ano}

We just saw in Section~\ref{subsec:Section-3-was-Ano} that various $(G,X)$-structures described in Section~\ref{sec:dyn-properties-hol} give rise (via the holonomy) to Anosov representations; these $(G,X)$-structures are of type~C or type~U (terminology of Section~\ref{subsec:intro-GXstruct}).
In this section, we study the converse direction.
Namely, given an Anosov representation $\rho : \Gamma\to G$, we~wish~to~find:
\begin{itemize}[leftmargin=0.8cm]
  \item homogeneous spaces $X=G/H$ on which $\Gamma$ acts properly discontinuously via~$\rho$; this will yield $(G,X)$-manifolds (or orbifolds) $M=\rho(\Gamma)\backslash X$ of type~C;
  \item proper open subsets $\mathcal{U}$ (\emph{domains of discontinuity}) of homogeneous spaces $X=G/H$ on which $\Gamma$ acts properly discontinuously via~$\rho$; this will yield $(G,X)$-manifolds (or orbifolds) $M=\rho(\Gamma)\backslash\mathcal{U}$ of type~U.
\end{itemize}
We discuss type~U in Sections \ref{subsec:DoD}--\ref{subsec:geom-Hitchin-max} and type~C in Section~\ref{subsec:complete-Anosov}.
One moti\-vation is to give a geometric meaning to the higher Teichm\"uller spaces of Section~\ref{subsec:Higher-Teich}.

\subsection{Cocompact domains of discontinuity} \label{subsec:DoD}

Domains of discontinuity with compact quotient have been constructed in several settings in the past ten years.

Barbot \cite{bar10} constructed such domains in the space $X$ of flags of~$\R^3$, for the Anosov representations to $G=\PSL(3,\R)$ of Example~\ref{ex:classical-Ano}.(3) and their small deformations.

Guichard--Wienhard \cite{gw12} developed a more general construction of cocompact domains of discontinuity in flag varieties $X$ for Anosov representations into semisimple Lie groups~$G$.
Here is one of their main results.
For $p\geq q\geq i\geq 1$, we denote~by $\mathcal{F}_i^{p,q}$ the closed subspace of the Grassmannian $\mathrm{Gr}_i(\R^{p+q})$ consisting of $i$-planes that are totally isotropic for the standard symmetric bilinear form $\langle\cdot,\cdot\rangle_{p,q}$ of signature~$(p,q)$.

\begin{theorem}[{Guichard--Wienhard \cite{gw12}}] \label{thm:gw-DoD}
Let $G=\PO(p,q)$ with $p\geq q$ and $X=\mathcal{F}_q^{p,q}$.
For any $P_1$-Anosov representation $\rho : \Gamma\to G\subset\PGL(\R^{p+q})$ with boundary map $\xi_1 : \partial_{\infty}\Gamma\to\mathcal{F}_1^{p,q}\subset\PP(\R^{p+q})$, the group $\rho(\Gamma)$ acts properly discontinuously with compact quotient on $\mathcal{U}_{\rho} := X \smallsetminus \mathscr{L}_{\rho}$, where
$$\mathscr{L}_{\rho} := \bigcup_{\eta\in\partial_{\infty}\Gamma} \{ W\in X=\mathcal{F}_q^{p,q} ~|~ \xi_1(\eta)\in W\}.$$
We have $\mathcal{U}_{\rho}\neq\emptyset$ as soon as $\dim(\partial_{\infty}\Gamma)<p-1$.
The homotopy type of $\rho(\Gamma)\backslash\mathcal{U}_{\rho}$ is constant as $\rho$ varies continuously among $P_1$-Anosov representations of $\Gamma$ in~$G$.
\end{theorem}

For $q=1$, we recover the familiar picture of Section~\ref{subsec:cc-rank1}: the set $\mathscr{L}_{\rho}$ is the limit set $\Lambda_{\rho(\Gamma)}\subset\partial_{\infty}\HH^p$, and $\mathcal{U}_{\rho}$ is the domain of discontinuity $\Omega_{\rho(\Gamma)}=\partial_{\infty}\HH^p\smallsetminus\Lambda_{\rho(\Gamma)}$.

Guichard--Wienhard \cite{gw12} used Theorem~\ref{thm:gw-DoD} to describe domains of discontinuity for various families of Anosov representations into other semisimple Lie groups~$G$.
Indeed, by \cite{gw12} an Anosov representation $\rho : \Gamma\to G$ can always be composed with a representation of $G$ into some $\PO(p,q)$ so as to become $P_1$-Anosov in $\PO(p,q)$.

Kapovich--Leeb--Porti \cite{klp15} developed a more systematic approach to the construction of domains of discontinuity in flag varieties.
They provided sufficient~conditions (expressed in terms of a notion of \emph{balanced} ideal in the Weyl group for the Bruhat order) on triples $(G,P,Q)$ consisting of a semisimple Lie group~$G$ and~two parabolic subgroups $P$ and~$Q$, so that $P$-Anosov representations into~$G$ admit cocompact domains of discontinuity in $G/Q$.
These domains are obtained by removing an explicit ``extended limit set'' $\mathscr{L}_{\rho}$ as in Theorem~\ref{thm:gw-DoD}.
The approach of Kapovich--Leeb--Porti is intrinsic: it does not rely on an embedding of $G$ into some~$\PO(p,q)$.

\subsection{Geometric structures for Hitchin and maximal representations} \label{subsec:geom-Hitchin-max}

Let $S$ be a closed orientable surface of genus $\geq 2$.
Recall (Example~\ref{ex:classical-Ano}) that Hitchin representations from $\Gamma=\pi_1(S)$ to $G=\PSL(d,\R)$ are $P_i$-Anosov for all $1\leq i\leq\nolinebreak d-\nolinebreak 1$; maximal representations from $\Gamma$ to $G=\PO(2,q)\subset\PGL(2+q,\R)$ are $P_1$-Anosov.

For $G=\PSL(2,\R)\simeq\PO(2,1)_0$, Hitchin representations and maximal representations of $\Gamma$ to~$G$ both coincide with the Fuchsian representations; they are the holonomy representations of hyperbolic structures on~$S$ (Example~\ref{ex:holonomy-Teich}).
In the setting of higher Teichm\"uller theory (see Section~\ref{subsec:Higher-Teich}), one could hope that Hitchin or maximal representations of $\Gamma$ to higher-rank~$G$ might also parametrize certain geometric structures on a manifold related to~$S$.
We saw in Section~\ref{subsec:Section-3-was-Ano} that this is indeed the case for $G=\PSL(3,\R)$: Hitchin represen\-tations of $\Gamma$ to $\PSL(3,\R)$ parametrize the convex projective structures~on~$S$~\cite{cg05,gol90}.
In an attempt to generalize this picture, we now outline constructions of domains of discontinuity for Hitchin representations to $\PSL(d,\R)$ with $d>3$, and maximal representations to $\PO(2,q)$ with $q>1$.
By classical considerations of cohomological dimension, such domains cannot be both cocompact and contractible; in Sections \ref{subsubsec:odd-Hitchin}--\ref{subsubsec:geom-max} below, we will prefer to forgo compactness to favor the nice geometry of convex domains.

\subsubsection{Hitchin representations for even $d=2n$} \label{subsubsec:even-Hitchin}

Let $(G,X)=(\PSL(2n,\R),\PP(\R^{2n}))$.
Hitchin representations to~$G$ do not preserve any properly convex open set $\Omega$ in~$X$, see \cite{dgk-proj-cc,zim}.
However, Guichard--Wienhard associated to them nonconvex $(G,X)$-structures on a closed manifold: by \cite{gw12}, if $\rho : \Gamma\to G$ is Hitchin with boundary map $\xi_n : \partial_{\infty}\Gamma\to\mathrm{Gr}_n(\R^{2n})$, then $X\smallsetminus \bigcup_{\eta\in\partial_{\infty}\Gamma} \xi_n(\eta)$ is a cocompact domain of discontinuity for~$\rho$, and the homotopy type of the quotient does not depend on~$\rho$.~So~far the topology and geometry of the quotient are understood only for $n=2$ (see \cite{wie-icm}).

\subsubsection{Hitchin representations for odd $d=2n+1$} \label{subsubsec:odd-Hitchin}

Hitchin representations to $G=\PSL(2n+1,\R)$ give rise to $(G,X)$-manifolds for at least two choices of~$X$.

One choice is to take $X$ to be the space of partial flags\,$(V_1\subset\nolinebreak V_{2n})$ of $\R^{2n+1}$ with $V_1$ a line and $V_{2n}$ a hyperplane: Guichard--Wienhard \cite{gw12} again constructed explicit cocompact domains of discontinuity in~$X$ in this setting.

Another choice is $X=\PP(\R^{2n+1})$: Hitchin representations in odd dimension are the holo\-nomies of convex projective manifolds, which are noncompact for $n>1$.

\begin{theorem}[\cite{dgk-proj-cc,zim}] \label{thm:odd-Hitchin-cc}
For any Hitchin representation $\rho : \Gamma\to\PSL(2n+1,\R)$, there is a $\rho(\Gamma)$-invariant properly convex open subset $\Omega$ of $\PP(\R^{2n+1})$ and a nonempty closed convex subset $\mathcal{C}$ of~$\Omega$ which has compact quotient by $\rho(\Gamma)$.
\end{theorem}

More precisely, if $\rho$ has boundary maps $\xi_1 : \partial_{\infty}\Gamma\to\mathrm{Gr}_1(\R^{2n+1})=\PP(\R^{2n+1})$ and $\xi_{2n} : \partial_{\infty}\Gamma\to\mathrm{Gr}_{2n}(\R^{2n+1})$, we may take $\Omega = \PP(\R^{2n+1}) \smallsetminus \bigcup_{\eta\in\partial_{\infty}\Gamma} \xi_{2n}(\eta)$ and $\mathcal{C}$ to be the convex hull of $\xi_1(\partial_{\infty}\Gamma)$ in~$\Omega$.
The group $\rho(\Gamma)$ acts properly discontinuously~on~$\Omega$ (Remark~\ref{rem:Hilbert}), and so $\rho(\Gamma)\backslash\Omega$ is a convex projective manifold, with a compact convex core $\rho(\Gamma)\backslash\CC$.
In other words, $\rho(\Gamma)$ is \emph{convex cocompact in $\PP(\R^{2n+1})$}, see Section~\ref{sec:proj-cc}.

\subsubsection{Maximal representations} \label{subsubsec:geom-max}

Maximal representations to $G=\PO(2,q)$ give rise to $(G,X)$-manifolds for at least two choices of~$X$.

One choice is $X=\mathcal{F}_2^{2,q}$ (also known as the space of \emph{photons} in the \emph{Einstein universe $\mathrm{Ein}^q$}): Theorem~\ref{thm:gw-DoD} provides cocompact domains of discontinuity for~$\rho$ in~$X$.
Collier--Tholozan--Toulisse \cite{ctt} recently studied the geometry of the associated quotient $(G,X)$-manifolds, and showed that they fiber over $S$ with~fiber~$\OO(q)/\OO(q-\nolinebreak 2)$.

Another choice is $X=\PP(\R^{2+q})$: by \cite{dgk-cc-Hpq}, maximal representations $\rho : \Gamma\to G$ are the holonomy representations of convex projective manifolds $\rho(\Gamma)\backslash\Omega$, which are noncompact for $q>1$ but still convex cocompact as in Section~\ref{subsubsec:odd-Hitchin}.
In fact $\Omega$ can be taken inside $\HH^{2,q-1}=\{[v]\in\PP(\R^{2+q})\,|\,\langle v,v\rangle_{2,q}<0\}$ (see \cite{dgk-cc-Hpq,ctt}), which is a pseudo-Riemannian analogue of the real hyperbolic space~in~signa\-ture $(2,q-1)$, and $\rho(\Gamma)$ is $\HH^{2,q-1}$-convex cocompact in the sense~of~Section~\ref{subsec:cc-Hpq}~below.

\subsection{Proper actions on full homogeneous spaces} \label{subsec:complete-Anosov}

In Sections \ref{subsec:DoD}--\ref{subsec:geom-Hitchin-max}, we mainly considered \emph{compact} spaces $X=G/H$ (flag varieties); these spaces cannot admit proper actions by infinite discrete groups, but we saw that sometimes they can contain domains of discontinuity $\mathcal{U}\subsetneq X$, yielding $(G,X)$-manifolds of type~U (terminology of Section~\ref{subsec:intro-GXstruct}).

We now consider \emph{noncompact} $X=G/H$.
Then Anosov representations $\rho :\nolinebreak \Gamma\to\nolinebreak G$ may give proper actions of~$\Gamma$ on the whole of $X=G/H$, yielding $(G,X)$-manifolds $\rho(\Gamma)\backslash X$ of type~C.
When $H$ is compact, this is not very interesting since all faithful and discrete representations to~$G$ give proper actions on~$X$.
However, when $H$ is noncompact, it may be remarkably difficult in general to find such representations giving proper actions on~$X$, which led to a rich literature (see \cite{ky05} and \cite[Intro]{kas-PhD}).

One construction for proper actions on~$X$ was initiated by Guichard--Wienhard~\cite{gw12} and developed further in~\cite{ggkw17bis}.
Starting from an Anosov representation $\rho : \Gamma\to G$, the idea is to embed $G$ into some larger semisimple Lie group~$G'$ so that $X=G/H$ identifies with a $G$-orbit in some flag variety $\mathcal{F}'$ of~$G'$, and then to find a cocompact domain of discontinuity $\mathcal{U}\supset X$ for $\rho$ in~$\mathcal{F}'$ by using a variant of Theorem~\ref{thm:gw-DoD}.
The action of $\rho(\Gamma)$ on~$X$ is then properly discontinuous, and $\rho(\Gamma)\backslash (\mathcal{U}\cap\overline{X})$ provides a compactification of $\rho(\Gamma)\backslash X$, which in many cases can be shown to be well-behaved.
Here is one application of this construction; see \cite{ggkw17bis} for other examples.

\begin{example}[\cite{ggkw17bis}] \label{ex:proper-higher-AdS}
Let $G=\PO(p,q)$ and $H=\OO(p,q-1)$ where $p>q\geq 1$.
For any $P_q$-Anosov representation $\rho :\nolinebreak\Gamma\to G\subset\PGL(p+q,\R)$, the group $\rho(\Gamma)$ acts properly discontinuously on $X = \HH^{p,q-1} = \{[v]\in\PP(\R^{p+q})\,|\,\langle v,v\rangle_{p,q}<0\} \simeq G/H$, and for torsion-free $\Gamma$ the complete $(G,X)$-manifold $\rho(\Gamma)\backslash X$ is topologically tame.
\end{example}

By \emph{topologically tame} we mean homeomorphic to the interior of a compact manifold with boundary.
For other compactifications of quotients of homogeneous spaces by Anosov representations, yielding topological tameness, see \cite{gkw,kl15,klp15}.

Another construction of complete $(G,X)$-manifolds for Anosov representations into reductive Lie groups~$G$ was given in \cite{ggkw17}, based on a properness criterion of Benoist \cite{ben96} and Kobayashi \cite{kob96}.
For simplicity we discuss it for $G=\PGL(d,\R)$.
As in Section~\ref{subsec:charact-Ano}, let $\mu_i(g)$ be the logarithm of the $i$-th singular value of a matrix $g\in\GL(d,\R)$; this defines a map $\mu = (\mu_1,\dots,\mu_d) : \PGL(d,\R)\to\R^d/\R(1,\dots,1)\simeq\R^{d-1}$.
The properness criterion of \cite{ben96,kob96} states that for two closed subgroups $H,\Gamma$ of $G=\PGL(d,\R)$, the action of $\Gamma$ on $G/H$ is properly discontinuous if and only if the set $\mu(\Gamma)$ ``drifts away at~infi\-nity from $\mu(H)$'', in the sense that for any $R>0$ we have $d_{\R^{d-1}}(\mu(\gamma),\mu(H))\geq R$ for all but finitely many $\gamma\in\Gamma$.
If $\Gamma$ is the image of an Anosov representation, then we can apply the implication (1)~$\Rightarrow$~(2) of Theorem~\ref{thm:charact-Ano} to see that the properness criterion is satisfied for many examples of~$H$.

\begin{example}[\cite{ggkw17}]
For $i=1$ (\resp $n$), the image of any $P_i$-Anosov representation to $G=\PSL(2n,\R)$ acts properly discontinuously on $X=G/H$ for $H=\SL(n,\C)$ (\resp $\SO(n+1,n-1)$).
\end{example}

\section{Convex cocompact projective structures} \label{sec:proj-cc}

In Sections \ref{sec:dyn-properties-hol} and~\ref{subsec:Section-3-was-Ano} we started from $(G,X)$-structures to produce Anosov representations, and in Section~\ref{sec:geom-for-Ano} we started from Anosov representations to produce $(G,X)$-structures.
We now discuss a situation, in the setting of convex projective geometry, in which the links between $(G,X)$-structures and Anosov representations are particularly tight and go in both directions, yielding a better understanding of both sides.
In Section~\ref{subsec:general-proj-cc} we will also encounter generalizations of Anosov representations, for finitely generated groups that are not necessarily Gromov hyperbolic.

\subsection{Convex cocompactness in higher real rank}

The results presented here are part of a quest to generalize the notion of rank-one convex cocompactness of Section~\ref{subsec:cc-rank1} to higher real rank.

The most natural generalization, in the setting of Riemannian symmetric spaces, turns out to be rather restrictive: Kleiner--Leeb~\cite{kl06} and Quint \cite{qui05} proved that if $G$ is a real simple Lie group of real rank $\geq 2$ and $K$ a maximal compact subgroup of~$G$, then any Zariski-dense discrete subgroup of~$G$ acting with compact quotient on some nonempty convex subset of $G/K$ is a uniform lattice in~$G$.

Meanwhile, we have seen in Section~\ref{subsec:ex-Anosov} that Anosov representations into higher-rank semisimple Lie groups $G$ have strong dynamical properties which nicely generalize those of rank-one convex cocompact representations (see \cite{klp-survey,gui-bourbaki}).
However, in general Anosov representations to~$G$ do not act with compact quotient on any nonempty convex subset of $G/K$, and it is not clear that Anosov representations should come with any geometric notion of convexity at all (see \eg Section~\ref{subsubsec:even-Hitchin}).

In this section, we shall see that Anosov representations in fact do come with convex structures.
We shall introduce several generalizations of convex cocompactness~to~hi\-gher real rank (which we glimpsed in Sections \ref{subsubsec:odd-Hitchin}--\ref{subsubsec:geom-max})~and~relate them to Anosov representations.
This is joint work with J.~Danciger and F.~Gu\'eritaud.

\subsection{Convex cocompactness in pseudo-Riemannian hyperbolic spaces} \label{subsec:cc-Hpq}

We start with a generalization of the hyperbolic quasi-Fuchsian manifolds of Example~\ref{ex:qF-intro} or the AdS quasi-Fuchsian manifolds of Section~\ref{subsec:AdS-qF}, where we replace the real hyperbolic space $\HH^3$ or its Lorentzian analogue $\AdS^3$ by their general pseudo-Riemannian analogue in signature $(p,q-1)$ for $p,q\geq 1$, namely
$$X = \HH^{p,q-1} = \big\{[v]\in\PP(\R^{p+q})\,|\,\langle v,v\rangle_{p,q}<0\big\}.$$
The symmetric bilinear form $\langle\cdot,\cdot\rangle_{p,q}$ of signature $(p,q)$ induces a pseudo-Riemannian structure of signature $(p,q-1)$ on~$X$, with isometry group $G=\PO(p,q)$ and constant negative sectional curvature (see \cite[\S\,2.1]{dgk-cc-Hpq}).
The following is not our original definition from \cite{dgk-cc-Hpq}, but an equivalent one from \cite[Th.\,1.25]{dgk-proj-cc}.

\begin{definition} \label{def:Hpq-cc}
A discrete subgroup $\Gamma$ of $G=\PO(p,q)$ is \emph{$\HH^{p,q-1}$-convex cocompact} if it preserves a properly convex open subset $\Omega$ of $X=\HH^{p,q-1}\subset\PP(\R^{p+q})$ and if it acts with compact quotient on some closed convex subset $\CC$ of~$\Omega$ with nonempty interior, whose ideal boundary $\partial_{\mathrm{i}}\CC := \overline{\CC} \smallsetminus \CC = \overline{\CC}\cap\partial X$ does not contain any nontrivial projective segment.
A representation $\rho : \Gamma\to G$ is \emph{$\HH^{p,q-1}$-convex cocompact} if its kernel is finite and its image is an $\HH^{p,q-1}$-convex cocompact subgroup of~$G$.
\end{definition}

Here $\overline{\CC}$ is the closure of $\CC$ in $\PP(\R^{p+q})$ and $\partial X$ the boundary of $X=\HH^{p,q-1}$ in $\PP(\R^{p+q})$.
For $\Gamma,\Omega,\mathcal{C}$ as in Definition~\ref{def:Hpq-cc}, the quotient $\Gamma\backslash\Omega$ is a $(G,X)$-manifold (or orbifold) (see Remark~\ref{rem:Hilbert}), which we shall call \emph{convex cocompact}; the subset $\Gamma\backslash\mathcal{C}$ is compact, convex, and contains all the topology, as in Sections \ref{subsec:cc-rank1} and~\ref{subsec:AdS-qF}.

There is a rich world of examples of convex cocompact $(G,X)$-manifolds, including direct generalizations of the quasi-Fuchsian manifolds of Sections \ref{subsec:cc-rank1} and~\ref{subsec:AdS-qF} (see \cite{bm12,dgk-cc-Hpq,dgk-proj-cc}) but also more exotic examples where the fundamental group is not necessarily realizable as a discrete subgroup of $\PO(p,1)$ (see \cite{dgk-cc-Hpq,lm}).

The following result provides links with Anosov representations.

\begin{theorem}[\cite{dgk-cc-Hpq,dgk-proj-cc}] \label{thm:main-Hpq}
For $p,q\geq 1$, let $\Gamma$ be an infinite discrete group and $\rho : \Gamma\to G=\PO(p,q)\subset\PGL(p+q,\R)$ a representation.
\begin{enumerate}
  \item\label{item:Hpq-cc-implies-Anosov} If $\rho$ is $\HH^{p,q-1}$-convex cocompact, then $\Gamma$ is Gromov hyperbolic and $\rho$ is $P_1$-Anosov.
  \item\label{item:Ano-convex-implies-Hpq-cc} Conversely, if $\Gamma$ is Gromov hyperbolic, if $\rho$ is $P_1$-Anosov, and if $\rho(\Gamma)$ preserves a properly convex open subset of $\PP(\R^{p+q})$, then $\rho$ is $\HH^{p,q-1}$-convex cocompact or $\HH^{q,p-1}$-convex cocompact.
  \item\label{item:Ano-connected-implies-Hpq-cc} If $\Gamma$ is Gromov hyperbolic with \emph{connected} boundary $\partial_{\infty}\Gamma$ and if $\rho$ is $P_1$-Anosov, then $\rho$ is $\HH^{p,q-1}$-convex cocompact or $\HH^{q,p-1}$-convex cocompact.
\end{enumerate}
\end{theorem}

In \eqref{item:Ano-convex-implies-Hpq-cc}--\eqref{item:Ano-connected-implies-Hpq-cc}, ``$\rho$ is $\HH^{q,p-1}$-convex cocompact'' is understood after identifying $\PO(p,q)$ with $\PO(q,p)$ and $\PP(\R^{p,q})\smallsetminus\overline{\HH^{p,q-1}}$ with $\HH^{q,p-1}$ under multiplication of $\langle\cdot,\cdot\rangle_{p,q}$~by~$-1$.
The\,case that\,$q=2$\,and\,$\Gamma$\,is\,a\,uniform\,lattice\,of\,$\PO(p,1)$\,is\,due\,to\,Barbot--M\'erigot\,\cite{bm12}.

The links between $\HH^{p,q-1}$-convex cocompactness and Anosov representations in Theorem~\ref{thm:main-Hpq} have several applications.

\smallskip

Applications to $(G,X)$-structures (see \cite{dgk-cc-Hpq,dgk-proj-cc}):\\
$\bullet$ $\HH^{p,q-1}$-convex cocompactness is stable under small deformations, because being Anosov is; thus the set of holonomy representations of convex cocompact $(G,X)$-structures on a given manifold~$M$ is open in $\Hom(\pi_1(M),G)$.\\
$\bullet$ Examples of convex cocompact $(G,X)$-manifolds can be obtained using classical families of Anosov representations: \eg Hitchin representations into $\PO(n+1,n)$ are $\HH^{n+1,n-1}$-convex cocompact for odd~$n$ and $\HH^{n,n}$-convex cocompact for even~$n$, and Hitchin representations into $\PO(n+1,n+1)$ are $\HH^{n+1,n}$-convex cocompact.
Maximal representations into $\PO(2,q)$ are $\HH^{2,q-1}$-convex cocompact, see Section~\ref{subsubsec:geom-max}.

\smallskip

Applications to Anosov representations:\\
$\bullet$ New examples of Anosov representations can be constructed from convex cocompact $(G,X)$-manifolds: \eg this approach is used in \cite{dgk-cc-Hpq} to prove that any Gromov hyperbolic right-angled Coxeter group in $d$ generators admits $P_1$-Anosov representations into $\PGL(d,\R)$.
This provides a large new class of hyperbolic groups admitting Anosov representations; these groups can have arbitrary large cohomological dimension and exotic boundaries (see \cite[\S\,1.7]{dgk-cc-Hpq} for references).
(Until now most known examples of Anosov representations were for surface groups or free~groups.)\\
$\bullet$ For $q=2$ and $\Gamma$ a uniform lattice of $\PO(p,1)_0$, Barbot \cite{bar15} used convex cocompact $(G,X)$-structures to prove that the connected component $\mathcal{T}$ of $\Hom(\Gamma,\PO(p,2))$ containing the natural inclusion $\Gamma\hookrightarrow\PO(p,1)_0\hookrightarrow\PO(p,2)$ consists entirely of Anosov representations.
This is interesting in the framework of Section~\ref{subsec:Higher-Teich}.

\subsection{Strong projective convex cocompactness} \label{subsec:strong-proj-cc}

We now consider a broader~notion of convex cocompactness, not involving any quadratic form.
Let $d\geq 2$.

\begin{definition} \label{def:strong-cc}
A discrete subgroup $\Gamma$ of $G=\PGL(d,\R)$ is \emph{strongly $\PP(\R^d)$-convex cocompact} if it preserves a strictly convex open subset $\Omega$ of $X=\PP(\R^d)$ with $C^1$ boundary and if it acts with compact quotient on some nonempty closed convex subset $\CC$ of~$\Omega$.
A representation $\rho : \Gamma\to G$ is \emph{strongly $\PP(\R^d)$-convex cocompact} if its kernel is finite and its image is a strongly $\PP(\R^d)$-convex cocompact subgroup of~$G$.
\end{definition}

The action of $\Gamma$ on~$\Omega$ in Definition~\ref{def:strong-cc} is a special case of a class of \emph{geometrically finite actions} introduced by~Crampon--Marquis~\cite{cm14}.
We use the adverb ``strongly'' to emphasize the strong regularity assumptions made on~$\Omega$.
In Definition~\ref{def:strong-cc} we say that the quotient $\Gamma\backslash\Omega$ is a \emph{strongly convex cocompact} projective manifold (or orbifold); the subset $\Gamma\backslash\mathcal{C}$ is again compact, convex, and contains all the topology.

Strongly $\PP(\R^d)$-convex cocompact representations include $\HH^{p,q-1}$-convex cocompact representations as in Section~\ref{subsec:cc-Hpq} (see \cite{dgk-cc-Hpq}), and the natural inclusion of groups dividing strictly convex open subsets of $\PP(\R^d)$ as in Section~\ref{subsec:conv-div}.
The following result generalizes Theorem~\ref{thm:main-Hpq}, and improves on earlier results of \cite{ben-cd-I,cm14}.

\begin{theorem}[\cite{dgk-proj-cc}] \label{thm:cc-Ano-PGL}
Let $\Gamma$ be an infinite discrete group and $\rho : \Gamma\to G=\PGL(d,\R)$ a representation such that $\rho(\Gamma)$ preserves a nonempty properly convex open subset of $X=\PP(\R^d)$.
Then $\rho$ is strongly $\PP(\R^d)$-convex cocompact if and only if $\Gamma$ is Gromov hyperbolic and $\rho$ is $P_1$-Anosov.
\end{theorem}

Another generalization of Theorem~\ref{thm:main-Hpq} was independently obtained by Zimmer \cite{zim}: it is similar to Theorem~\ref{thm:cc-Ano-PGL}, but involves a slightly different notion of convex cocompactness and assumes $\rho(\Gamma)$ to act irreducibly on $\PP(\R^d)$.

\smallskip

Applications of Theorem~\ref{thm:cc-Ano-PGL} include:\\
$\bullet$ Examples of strongly convex cocompact projective manifolds using classical Anosov representations (\eg Hitchin representations into $\PSL(2n+1,\R)$ as in Section~\ref{subsubsec:odd-Hitchin}).\\
$\bullet$ In certain cases, a better understanding of the set of Anosov representations of a Gromov hyperbolic group~$\Gamma$ inside a given connected component of $\Hom(\Gamma,G)$: \eg for an irreducible hyperbolic right-angled Coxeter group~$\Gamma$, it is proved in \cite{dgk-racg-cc}, using Theorem~\ref{thm:cc-Ano-PGL} and Vinberg's theory \cite{vin71}, that $P_1$-Anosov representations form the full interior of the space of faithful and discrete representations of $\Gamma$ as a reflection group in $G=\PGL(d,\R)$.

\smallskip

For a Gromov hyperbolic group $\Gamma$ and a $P_1$-Anosov representation $\rho : \Gamma\to G=\PGL(d,\R)$, the group $\rho(\Gamma)$ does not always preserve a properly convex open subset of $X=\PP(\R^d)$: see Section~\ref{subsubsec:even-Hitchin}.
However, as observed in \cite{zim}, $\rho$ can always be composed with the embedding $\iota : G\hookrightarrow\PGL(V)$ described in Example~\ref{ex:conv-div}.(2), for $V=\mathrm{Sym}(d,\R)\simeq\R^{d(d+1)}$; then $\iota\circ\rho(\Gamma)$ preserves a properly convex open subset in $\PP(V)$.
The composition $\iota\circ\rho : \Gamma\to\PGL(V)$ is still $P_1$-Anosov by \cite{gw12}, and it is strongly $\PP(V)$-convex cocompact by Theorem~\ref{thm:cc-Ano-PGL}.
More generally, using \cite{gw12}, any Anosov representation to any semisimple Lie group can always be composed with an embedding into some $\PGL(V)$ so as to become strongly $\PP(V)$-convex cocompact.

\subsection{Projective convex cocompactness in general} \label{subsec:general-proj-cc}

We now introduce an even broader notion of convex cocompactness, where we remove the strong regularity assumptions on~$\Omega$ in Definition~\ref{def:strong-cc}.
This yields a large class of convex projective manifolds, whose fundamental groups are not necessarily Gromov hyperbolic.
Their holonomy representations are generalizations of Anosov representations, sharing some of their desirable properties (Theorem~\ref{thm:general-cc}).
This shows that Anosov representations are not the only way to successfully generalize rank-one convex cocompactness to higher real rank.

\begin{definition}[\cite{dgk-proj-cc}] \label{def:cc-general}
A discrete subgroup $\Gamma$ of $G=\PGL(d,\R)$ is \emph{$\PP(\R^d)$-convex cocompact} if it preserves a properly convex open subset $\Omega$ of $X=\PP(\R^d)$ and if it acts with compact quotient on some ``large enough'' closed convex subset $\CC$ of~$\Omega$.
A representation $\rho : \Gamma\to G$ is \emph{$\PP(\R^d)$-convex cocompact} if its kernel is finite and its image is a $\PP(\R^d)$-convex cocompact subgroup of~$G$.
\end{definition}

In Definition~\ref{def:cc-general}, by ``$\CC$ large enough'' we mean that all accumulation points of all $\Gamma$-orbits of~$\Omega$ are contained in the boundary of $\CC$ in $X=\PP(\R^d)$.
If we did not impose this (even if we asked $\CC$ to have nonempty interior), then the notion of $\PP(\R^d)$-convex cocompactness would not be stable under small deformations:~see~\cite{dgk-proj-cc,dgk-cc-ex}.
In Definition~\ref{def:cc-general} we call $\Gamma\backslash\Omega$ a \emph{convex cocompact} projective manifold (or orbifold).

The class of $\PP(\R^d)$-convex cocompact representations includes all strongly $\PP(\R^d)$-convex cocompact representations as in Section~\ref{subsec:strong-proj-cc}, hence all $\HH^{p,q-1}$-convex cocompact representations as in Section~\ref{subsec:cc-Hpq}.
In fact, the following holds.

\begin{proposition}[\cite{dgk-proj-cc}] \label{prop:cc-strong-cc}
Let $\Gamma$ be an infinite discrete group.
A representation $\rho : \Gamma\to G=\PGL(d,\R)$ is strongly $\PP(\R^d)$-convex cocompact (Definition~\ref{def:strong-cc}) if and only if it is $\PP(\R^d)$-convex cocompact (Definition~\ref{def:cc-general}) and $\Gamma$ is Gromov hyperbolic.
\end{proposition}

This generalizes a result of Benoist \cite{ben-cd-I} on divisible convex sets.
Together with Theorem~\ref{thm:cc-Ano-PGL}, Proposition~\ref{prop:cc-strong-cc} shows that $\PP(\R^d)$-convex cocompact representations are generalizations of Anosov representations, to a larger class of finitely generated groups $\Gamma$ which are not necessarily Gromov hyperbolic.
These representations still enjoy the following good properties.

\begin{theorem}[\cite{dgk-proj-cc}] \label{thm:general-cc}
Let $\Gamma$ be an infinite discrete group and $\rho : \Gamma\to G=\PGL(d,\R)$ a $\PP(\R^d)$-convex cocompact representation.
Then
\begin{enumerate}
  \item $\rho$ is a quasi-isometric embedding;
  \item\label{item:cc-stable} there is a neighborhood of~$\rho$ in $\Hom(\Gamma,G)$ consisting entirely of faithful and discrete $\PP(\R^d)$-convex cocompact representations;
  \item $\rho$ is $\PP((\R^d)^*)$-convex cocompact;
  \item\label{item:cc-include} $\rho$ induces a $\PP(\R^D)$-convex cocompact representation for any $D\geq d$ (by lifting $\rho$ to a representation to $\SL^{\pm}(d,\R)$ and composing it with the natural inclusion $\SL^{\pm}(d,\R)\hookrightarrow\SL^{\pm}(D,\R)$).
\end{enumerate}
\end{theorem}

In order to prove \eqref{item:cc-stable}, we show that the representations of Theorem~\ref{thm:general-cc} are exactly the holonomy representations of compact convex projective manifolds with strictly convex boundary \cite{dgk-proj-cc}; we can then apply the deformation theory of~\cite{clt}.

Groups that are $\PP(\R^d)$-convex cocompact but not strongly $\PP(\R^d)$-convex cocompact include all groups dividing a properly convex, but not strictly convex, open subset of $X=\PP(\R^d)$ as in Section~\ref{subsec:conv-div}, as well as their small deformations in $\PGL(D,\R)$ for $D\geq d$ (Theorem~\ref{thm:general-cc}.\eqref{item:cc-stable}--\eqref{item:cc-include}).
Such nontrivial deformations exist: \eg for $d=4$ we can always bend along tori or Klein bottle subgroups \cite{ben-cd-IV}.
There seems to be a rich world of examples beyond this, which is just starting to be unveiled, see \cite{dgk-proj-cc,dgk-racg-cc,dgk-cc-ex}.
It would be interesting to understand the precise nature of the corresponding abstract groups~$\Gamma$, and how the dynamics of $\PP(\R^d)$-convex cocompact representations generalize that of Anosov representations.

\section{Complete affine structures} \label{sec:hol-are-Ano}

In Sections \ref{sec:dyn-properties-hol} to~\ref{sec:proj-cc} we always considered semisimple, or more generally reductive, Lie groups~$G$.
We now discuss links between $(G,X)$-structures and representations of discrete groups into~$G$ in an important case where $G$ is not reductive: namely~$G$ is the group $\mathrm{Aff}(\R^d)=\GL(d,\R)\ltimes\nolinebreak\R^d$ of invertible affine transformations of $X=\nolinebreak\R^d$.
We shall see in Section~\ref{subsec:complete-affine-crit} that for $d=3$ the holonomy representations of certain complete (\ie type~C in Section~\ref{subsec:intro-GXstruct}) $(G,X)$-structures are characterized by a uniform contraction condition, which is also an \emph{affine Anosov} condition; we shall briefly mention partial extensions to $d>3$, which are currently being~explored.

\subsection{Brief overview: understanding complete affine manifolds} \label{subsec:overview-complete-affine}

Let $(G,X)=(\mathrm{Aff}(\R^d),\R^d)$.
This section is centered around \emph{complete affine manifolds}, \ie $(G,X)$-manifolds of the form $M=\rho(\Gamma)\backslash X$ where $\Gamma\simeq\pi_1(M)$ is a discrete group and $\rho : \Gamma\to G$ a faithful representation through which $\Gamma$ acts properly discontinuously and freely on~$X=\R^d$.
The study of such representations has a rich history through the interpretation of their images as \emph{affine crystallographic groups}, \ie symmetry groups of affine tilings of~$\R^d$, possibly with noncompact tiles; see \cite{abe01} for a detailed survey.
The compact and noncompact cases are quite different.

For a compact complete affine manifold $M$, Auslander \cite{aus64} conjectured that~$\pi_1(M)$ must always be virtually (\ie up to finite index) polycyclic.
This extends a classical theorem of Bieberbach on~affi\-ne Euclidean isometries.
The conjecture is proved for $d\leq 6$ \cite{ams12,fg83}, but remains widely open for $d\geq 7$, despite partial results~(see~\cite{abe01}).

In contrast, in answer to a question of Milnor \cite{mar84}, there exist noncompact complete affine manifolds $M$ for which $\pi_1(M)$ is \emph{not} virtually polycyclic.
The first examples were constructed by Margulis~\cite{mar84} for $d=3$, with $\pi_1(M)$ a nonabelian free group.
In these examples the holonomy representation takes values in $\OO(2,1)\ltimes\R^3$ (this is always the case when $\pi_1(M)$ is not virtually polycyclic \cite{fg83}), hence $M$ inherits a flat Lorentzian structure.
Such manifolds are called \emph{Margulis spacetimes}.
They have a rich geometry and have been much studied since the 1990s, most prominently by Charette, Drumm, Goldman, Labourie, and Margulis.
In particular, the questions of the topological tameness of Margulis spacetimes and of the existence of nice fundamental domains in $X=\R^3$ (bounded by piecewise linear objects called \emph{crooked planes}) have received much attention: see \eg \cite{cdg15,cg17,dgk16,dgk16bis,dru92,dg99}.
See also \cite{ams12,glm09,smi16} for higher-dimensional analogues $M$ with $\pi_1(M)$ a free group.

Following \cite{dgk16,dgk16bis} (see \cite{sch-bourbaki}), a convenient point of view for understanding Margulis spacetimes is to regard them as ``infinitesimal analogues'' of complete AdS~ma\-nifolds.
In order to describe this point of view, we first briefly discuss the AdS~case.

\subsection{Complete AdS manifolds} \label{subsec:complete-AdS}

As in Section~\ref{subsec:AdS-qF}, let $(G,X)=(\PO(2,2),\AdS^3)$, and view $X$ as the group $\underline{G}=\PSL(2,\R)$ and the identity component $G_0$ of~$G$ as $\underline{G}\times\underline{G}$ acting on $X\simeq\underline{G}$ by right and left multiplication.
We consider $(G,X)$-manifolds of the form $M=\rho(\Gamma)\backslash X$ where $\Gamma\simeq\pi_1(M)$ is an infinite discrete group and $\rho = (\rho_L,\rho_R) : \Gamma\to\underline{G}\times\underline{G}\subset G$ a faithful representation through which $\Gamma$ acts properly discontinuously and freely on~$X$.
Not all faithful and discrete $\rho=(\rho_L,\rho_R)$ yield properly discontinuous actions on~$X$: \eg if $\rho_L=\rho_R$, then $\rho$ has a global fixed point, precluding properness.
However, the following properness criteria hold.
We denote by $\lambda(g):=\inf_{x\in\HH^2} d_{\HH^2}(x,g\cdot x)\geq 0$ the translation length of $g\in\underline{G}$ in~$\HH^2$.

\begin{theorem}[\cite{kas-PhD,ggkw17}] \label{thm:prop-crit-AdS3}
Let $G=\PO(2,2)$ and $\underline{G}=\PSL(2,\R)$.
Consider~a~discrete group $\Gamma$ and a representation $\rho = (\rho_L,\rho_R) : \Gamma\to\underline{G}\times\underline{G}\subset G$ with $\rho_L$ convex~co\-compact.
The following are equivalent, up to switching $\rho_L$ and~$\rho_R$~in~both~\eqref{item:contract-lambda}~and~\eqref{item:contract-Lip}:
\begin{enumerate}
  \item\label{item:proper-AdS} the action of $\Gamma$ on $X=\AdS^3\simeq\underline{G}$ via~$\rho$ is properly discontinuous;
  \item\label{item:contract-lambda} there exists $C<1$ such that $\lambda(\rho_R(\gamma))\leq C \lambda(\rho_L(\gamma))$ for all $\gamma\in\Gamma$;
  \item\label{item:contract-Lip} there is a $(\rho_L,\rho_R)$-equivariant Lipschitz map $f : \HH^2\to\HH^2$ with $\mathrm{Lip}(f)<1$;
  \item\label{item:GxG-PSL2-Ano} $\Gamma$ is Gromov hyperbolic and $\rho : \Gamma\to G\subset\PGL(4,\R)$ is $P_2$-Anosov.
\end{enumerate}
\end{theorem}

The equivalences \eqref{item:proper-AdS}~$\Leftrightarrow$~\eqref{item:contract-lambda}~$\Leftrightarrow$~\eqref{item:contract-Lip}, proved in \cite{kas-PhD}, have been generalized in \cite{gk17} to $\underline{G}=\PO(n,1)$ for any $n\geq 2$, allowing $\rho_L$ to be geometrically finite instead of convex cocompact.
These equivalences state that $\rho=(\rho_L,\rho_R)$ acts properly discontinuously on $X=\AdS^3\simeq\underline{G}$ if and only if, up to switching the two factors, $\rho_L$ is faithful and discrete and $\rho_R$ is ``uniformly contracting'' with respect to~$\rho_L$.
The equivariant map $f$ in~\eqref{item:contract-Lip} provides an explicit fibration in circles of $\rho(\Gamma)\backslash X$ over the hyperbolic surface $\rho_L(\Gamma)\backslash\HH^2$, see \cite{gk17}.
We refer to \cite{dgk-cox,dt16,gk17,gkw15,ll17,sal00} for many examples, to \cite{tho17} for a classification in the compact AdS case, and to \cite{gk17,kas-PhD} for links with Thurston's asymmetric metric on Teichm\"uller space \cite{thu86}.

The equivalences \eqref{item:proper-AdS}~$\Leftrightarrow$~\eqref{item:contract-lambda}~$\Leftrightarrow$~\eqref{item:GxG-PSL2-Ano}, proved in \cite{ggkw17}, generalize to $\underline{G}=\PO(n,1)$, $\mathrm{PU}(n,1)$, or $\mathrm{Sp}(n,1)$; the Anosov condition is then expressed in $\PGL(2n+2,\mathbb{K})$ where $\mathbb{K}$ is $\R$, $\C$, or the quaternions.
As an application \cite{ggkw17}, the set of holonomy representations of \emph{complete} $(\underline{G}\!\times\!\underline{G},\,\underline{G})$-structures on a compact manifold $M$ is open in the set of holonomy representations of all possible $(\underline{G}\!\times\!\underline{G},\,\underline{G})$-structures on~$M$.
By \cite{tho15}, it is also closed, which gives evidence for an open conjecture stating that all $(\underline{G}\!\times\!\underline{G},\,\underline{G})$-structures on~$M$ should be complete (\ie obtained as quotients of~$\widetilde{\underline{G}}$).

\subsection{Complete affine manifolds} \label{subsec:complete-affine-crit}

We now go back to $(G,X)=(\mathrm{Aff}(\R^d),\R^d)$, loo\-king for characterizations of holonomy representations of complete affine manifolds, \ie representations into~$G$ yielding properly discontinuous actions on~$X$.

We first note that any representation from a group $\Gamma$ to the nonreductive Lie group $G=\GL(d,\R)\ltimes\R^d$ is of the form $\rho=(\rho_L,u)$ where $\rho_L : \Gamma\to\GL(d,\R)$ (linear part) is a representation to $\GL(d,\R)$ and $u : \Gamma\to\R^d$ (translational part) a $\rho_L$-cocycle, meaning $u(\gamma_1\gamma_2) = u(\gamma_1) + \rho_L(\gamma_1)\cdot u(\gamma_2)$ for all $\gamma_1,\gamma_2\in\Gamma$.

We focus on the case $d=3$ and $\rho_L$ with values in $\OO(2,1)$.
Let us briefly indicate how, following \cite{dgk16,dgk16bis}, the Margulis spacetimes of Section~\ref{subsec:overview-complete-affine} are ``infinitesimal versions'' of the complete AdS manifolds of Section~\ref{subsec:complete-AdS}.
Let $\underline{G}=\OO(2,1)_0\simeq\PSL(2,\R)$ be the group of isometries of~$\HH^2$.
Its Lie algebra $\smash{\underline{\g}}\simeq\R^3$ is the set of ``infinitesimal isometries'' of~$\HH^2$, \ie Killing vector fields on~$\HH^2$.
Here are some properness criteria.

\begin{theorem}[\cite{glm09,dgk16}] \label{thm:prop-crit-R21}
Let $G=\mathrm{Aff}(\R^3)$ and $\underline{G}=\OO(2,1)_0\simeq\PSL(2,\R)$.~Consider a discrete group $\Gamma$ and a representation $\rho = (\rho_L,u) : \Gamma\to\underline{G}\ltimes\smash{\underline{\g}}\subset G$~with~$\rho_L$~con\-vex cocompact.
The following are equivalent, up to replacing $u$ by $-u$ in both~\eqref{item:contract-inf-lambda}~and~\eqref{item:contract-inf-Lip}:
\begin{enumerate}
  \item\label{item:proper-R21} the action of $\Gamma$ via $\rho=(\rho_L,u)$ on $X=\R^3\simeq\smash{\underline{\g}}$ is properly discontinuous;
  \item\label{item:contract-inf-lambda} there exists $c<0$ such that $\frac{\mathrm{d}}{\mathrm{d}t}|_{t=0}\,\lambda(e^{u(\gamma)}\rho_L(\gamma))\leq c\,\lambda(\rho_L(\gamma))$ for all $\gamma\in\Gamma$;
  \item\label{item:contract-inf-Lip} there is a $(\rho_L,u)$-equivariant vector field $Y$ on~$\HH^2$ with\,``lipschitz''\,constant\,$<\nolinebreak 0$.
\end{enumerate}
\end{theorem}

The equivalence \eqref{item:proper-R21}~$\Leftrightarrow$~\eqref{item:contract-inf-lambda} is a reinterpretation, based on \cite{gm00}, of a celebrated result of Goldman--Labourie--Margulis \cite{glm09}.
The equivalence \eqref{item:proper-R21}~$\Leftrightarrow$~\eqref{item:contract-inf-Lip} is proved~in~\cite{dgk16}.

These equivalences are ``infinitesimal versions'' of the equivalences \eqref{item:proper-AdS}~$\Leftrightarrow$~\eqref{item:contract-lambda}~$\Leftrightarrow$~\eqref{item:contract-Lip} of Theorem~\ref{thm:prop-crit-AdS3}.
Indeed, as explained in \cite{dgk16}, we can see the $\rho_L$-cocycle $u : \Gamma\to\underline{\g}$ as an  ``infinitesimal deformation'' of the holonomy representation $\rho_L$ of the hyperbolic surface (or orbifold) $S=\rho_L(\Gamma)\backslash\HH^2$; Condition~\eqref{item:contract-inf-lambda} states that closed geodesics on~$S$ get uniformly shorter under this infinitesimal deformation.
We can see a $(\rho_L,u)$-equivariant vector field $Y$ on~$\HH^2$ as an ``infinitesimal deformation'' of the developing map of the hyperbolic surface~$S$; Condition~\eqref{item:contract-inf-Lip} states that any two points of~$\HH^2$ get uniformly closer compared to their initial distance.
Thus Theorem~\ref{thm:prop-crit-R21} states that $\rho=(\rho_L,u)$ acts properly discontinuously on $X=\R^3\simeq\smash{\underline{\g}}$ if and only if the infinitesimal deformation~$u$, up to replacing it by $-u$, is ``uniformly contracting''.

The vector field $Y$ in~\eqref{item:contract-inf-Lip} provides an explicit fibration in lines of the Margulis spacetime $\rho(\Gamma)\backslash X$ over the hyperbolic surface~$S$, and this can be used to define a \emph{geometric transition} from complete AdS manifolds to Margulis spacetimes,~see~\cite{dgk16}.

In Theorem~\ref{thm:prop-crit-AdS3}, the ``uniform contraction'' characterizing properness was in fact an Anosov condition, encoding strong dynamics on a certain flag variety.
It is natural to expect that something similar should hold in the setting of Theorem~\ref{thm:prop-crit-R21}.
For this, a notion of \emph{affine Anosov representation} into $\OO(2,1)\ltimes\R^3$ was recently introduced by Ghosh \cite{gho15} and extended to $\OO(n+1,n)\ltimes\R^d\subset\mathrm{Aff}(\R^d)=G$ for~any $d=2n+1\geq 3$ by Ghosh--Treib \cite{gt17}; the definition is somewhat analogous~to~Section~\ref{subsec:def-Ano} but uses affine bundles and their sections.
By \cite{gho15,gt17}, given a $P_n$-Anosov representation $\rho_L : \Gamma\to\OO(n+1,n)$ and a cocycle $u : \Gamma\to\R^d$, the action of $\Gamma$ on $X=\R^d$ via $\rho=(\rho_L,u)$ is properly discontinuous if and only if $\rho$ is affine Anosov.

Theorem~\ref{thm:prop-crit-R21} was recently generalized in \cite{dgk-cox} as follows: for $\underline{G}=\OO(p,q)$ with $p,q\geq\nolinebreak 1$, consider a discrete group~$\Gamma$, a faithful and discrete representation $\rho_L : \Gamma\to\underline{G}$, and a $\rho_L$-cocycle $u : \Gamma\to\smash{\underline{\g}}$; then the action of $\Gamma$ on~$\smash{\underline{\g}}$ via $\rho = (\rho_L,u) : \Gamma\to\mathrm{Aff}(\smash{\underline{\g}})$ is properly discontinuous as soon as $u$ satisfies a uniform contraction property in the pseudo-Riemannian hyperbolic space $\HH^{p,q-1}$ of Section~\ref{subsec:cc-Hpq}.
This allowed for the construction in \cite{dgk-cox} of the first examples of irreducible complete affine manifolds $M$ such that $\pi_1(M)$ is neither virtually polycyclic nor virtually free: $\pi_1(M)$ can in fact be any right-angled Coxeter group.
It would be interesting to understand the links with a notion of affine Anosov representation in this setting.

\section{Concluding remarks} \label{sec:conclusion}

By investigating the links between the geometry of $(G,X)$-structures on manifolds and the dynamics of their holonomy representations, we have discussed only a small part of a very active area of research.

We have described partial answers to Problem~\ref{problemA} for several types of model geometries $(G,X)$.
However, Problem~\ref{problemA} is still widely open in many contexts.
As an illustration, let us mention two major open conjectures on closed affine manifolds (in addition to the Auslander conjecture of Section~\ref{subsec:overview-complete-affine}): the Chern conjecture states that if a closed $d$-manifold~$M$ admits an $(\mathrm{Aff}(\R^d),\R^d)$-structure, then its Euler characteristic must be zero; the Markus conjecture states that an $(\mathrm{Aff}(\R^d),\R^d)$-structure on~$M$ is complete if and only if its holonomy representation takes values in $\SL(d,\R)\ltimes\R^d$.
See \cite{kli17} and references therein for recent progress on this.

We have seen that Anosov representations from Gromov hyperbolic groups to semisimple Lie groups provide a large class of representations answering Problem~\ref{problemB}.
However, not much is known beyond them.
One further class, generalizing Anosov representations to finitely generated groups $\Gamma$ which are not necessarily Gromov hyperbolic, is the class of $\PP(\R^d)$-convex cocompact representations into $\PGL(d,\R)$ of Section~\ref{subsec:general-proj-cc}; it would be interesting to understand this class better in the framework of Problem~\ref{problemB}, see Section~\ref{subsec:general-proj-cc} and \cite[Appendix]{dgk-proj-cc}.
As another generalization of Anosov representations, it is natural to look for a class of representations of relatively hyperbolic groups into higher-rank semisimple Lie groups which would bear similarities to geometrically finite representations into rank-one groups, with cusps allowed: see \cite[\S\,5]{kl17} for a conjectural picture.
Partial work in this direction has been done in the convex projective setting, see \cite{cm14}.

To conclude, here are two open questions which we find particularly interesting.

\subsection*{Structural stability}

Sullivan \cite{sul85} proved that a structurally stable, nonrigid subgroup of $G=\PSL(2,\C)$ is always Gromov hyperbolic and convex cocompact in~$G$.
It is natural to ask if this may be extended to subgroups of higher-rank semisimple Lie groups $G$ such as $\PGL(d,\R)$ for $d\geq 3$, for instance with ``convex cocompact'' replaced by ``Anosov''.
In Section~\ref{subsec:general-proj-cc} we saw that there exist nonrigid, structurally stable~sub\-groups of $G=\PGL(d,\R)$ which are \emph{not} Gromov hyperbolic, namely groups that are $\PP(\R^d)$-convex cocompact but not strongly $\PP(\R^d)$-convex cocompact (Definitions \ref{def:strong-cc} and~\ref{def:cc-general}).
However, does a Gromov hyperbolic, nonrigid, structurally stable, discrete subgroup of~$G$ always satisfy some Anosov property?

\subsection*{Abstract groups admitting Anosov representations} \label{subsec:groups-admitting-Ano}

Which linear hyperbolic groups admit Anosov representations to some semisimple Lie group?
Classical examples~include surface groups, free groups, and more generally rank-one convex cocompact groups, see Section~\ref{subsec:ex-Anosov}.
By \cite{dgk-cc-Hpq}, all Gromov hyperbolic right-angled Coxeter groups (and all groups commensurable to them) admit Anosov representations, see Section~\ref{subsec:cc-Hpq}.
On the other hand, if a hyperbolic group admits an Anosov representation, then its Gromov flow (see Section~\ref{subsec:def-Ano}) must satisfy strong dynamical properties, which may provide an obstruction: see the final remark of \cite[\S\,1]{bcls15}.
It would be interesting to have further concrete examples of groups admitting or not admitting Anosov representations.


\end{document}